\documentclass[aoas,preprint]{imsart}

\RequirePackage[OT1]{fontenc}
\RequirePackage{amsthm,amsmath,amssymb}
\RequirePackage{natbib}
\RequirePackage[colorlinks,citecolor=blue,urlcolor=blue]{hyperref}
\usepackage{natbib}
\usepackage{graphicx}
\usepackage{epsfig}
\usepackage{psfrag}
\usepackage{wrapfig}
\usepackage[all]{xy}
\startlocaldefs
\newtheorem{lemma}{Lemma}[section]
\newtheorem{prop}{Proposition}[section]
\newtheorem{corollary}{Corollary}[section]
\newtheorem{cond}{Condition}
\def\Var{{\rm Var}\,}
\def\E{{\rm E}\,}
\def\arg{{\rm arg}\,}
\def\Cov{{\rm Cov}\,}
\def\N{{\rm N}\,}
\newcommand{\Z}{\mathbf{Z}}
\newcommand{\I}{\mathbf{I}}

\newcommand{\z}{\mathbf{z}}

\newcommand{\bpi}{\boldsymbol{\pi}}
\newcommand{\arrowp}{\stackrel{p}{\rightarrow}}

\newcommand{\bP}{\mathbf{P}}

\newcommand\independent{\protect\mathpalette{\protect\independenT}{\perp}}
\def\independenT#1#2{\mathrel{\rlap{$#1#2$}\mkern2mu{#1#2}}}
\usepackage{graphicx}
\endlocaldefs
\begin{document}

\sloppy

\begin{frontmatter}
\title{Estimating Average Causal Effects Under General Interference, with Application to a Social Network Experiment}
\runtitle{Causal Effects Under Interference}
%
\begin{aug}
\author{\fnms{Peter M.} \snm{Aronow}\thanksref{m1}\ead[label=e1]{peter.aronow@yale.edu}}
\and
\author{\fnms{Cyrus} \snm{Samii}\thanksref{m3}\ead[label=e2]{cds2083@nyu.edu}}

\runauthor{Aronow and Samii}
\affiliation{Yale University\thanksmark{m1}
and New York University\thanksmark{m3}}
\address{Peter M. Aronow \\
Departments of Political Science and Biostatistics \\ 
Yale University \\  77 Prospect Street \\ New Haven, CT 06520 \\ United States \\
\printead{e1}}
\address{Cyrus Samii \\
Politics Department\\
New York University\\
19 West 4th Street \\ 
 New York, NY 10012\\ 
 United States \\
\printead{e2}\\}
\end{aug}
%
%
\begin{abstract}
This paper presents a randomization-based framework for estimating causal effects under interference between units, motivated by challenges that arise in analyzing experiments on social networks.  The framework integrates three components: (i) an experimental design that defines the probability distribution of treatment assignments, (ii) a mapping that relates experimental treatment assignments to exposures received by units in the experiment, and (iii) estimands that make use of the experiment to answer questions of substantive interest.  We develop the case of estimating average unit-level causal effects from a randomized experiment with interference of arbitrary but known form. The resulting estimators are based on inverse probability weighting. We provide randomization-based variance estimators that account for the complex clustering that can occur when interference is present.  We also establish consistency and asymptotic normality under local dependence assumptions.  We discuss refinements including covariate-adjusted effect estimators and ratio estimation. We evaluate empirical performance in realistic settings with a naturalistic simulation using social network data from American  schools. We then present results from a field experiment on the spread of anti-conflict norms and behavior among school students.  
\end{abstract}



\end{frontmatter}

\section{Introduction} \label{interfereintro}

We develop methods for analyzing an experiment in which treatments are applied to individuals in a social network and causal effects are hypothesized to transmit to peers through the network.  Experimental and observational studies often involve treatments with effects that ``interfere'' \citep{cox58}  across units through spillover or other forms of dependency. Such interference is sometimes considered a nuisance, and researchers may strive to design studies that isolate units as much as possible from interference. However, such designs are not always possible. Furthermore, researchers may be interested in estimation of the spillover effects themselves, as these effects may be of substantive importance.  Other applications share structural similarities to the social network case. For example, an urban renewal program applied to one town may divert capital from other towns, in which case the overall effect of the program may be ambiguous. Treatment effects may carry over from one time period to another, and units have some chance of receiving treatment at any one of a set of points in time. In these cases, we need methods to estimate effects of both direct and indirect exposure to a treatment.   Moreover, researchers may be interested in understanding how such indirect effects vary depending on individuals' characteristics.

This paper 
presents a 
general, randomization-based 
framework
for estimating 
causal effects under these and other forms of interference. 
Interference represents a departure from the traditional 
causal inference scenario
wherein units are 
assigned directly to treatment or control, and the potential outcomes that would be observed for a unit in either the treatment or control condition are fixed \citep{colefrangakis2009} and do not depend on the overall set of treatment assignments.  The latter condition is what \citet{rubin1990} refers to as the ``stable unit treatment value assumption'' (SUTVA).  In the examples above, the traditional 
scenario
is clearly 
an
inadequate
characterization, 
as SUTVA would be violated.  
A more sophisticated 
characterization
of treatment exposure 
and associated potential outcomes
must be specified.
For the school field experiment, program participation was randomly assigned, but encouragement to support anti-conflict norms could come from direct participation in the program as well as indirect exposure via social network peers that participated in the program.

We start with theoretical results for  
an
estimation framework 
that
consists of three components: (i) the experimental (or quasi-experimental) ``design,'' which characterizes precisely the probability distribution of treatments {\it assigned}; (ii) an ``exposure mapping,'' which relates treatments assigned to exposures {\it received}; and (iii) a set of causal estimands selected to make use of the experiment to answer questions of substantive interest.  
For the case of a randomized experiment under arbitrary but known forms of interference, we provide 
unbiased estimators of average unit-level causal effects induced by treatment exposure.  We also provide estimators for the randomization variance of the estimated average causal effects.  These variance estimators are assured of being conservative (that is, nonnegatively biased). We 
establish conditions for consistency and 
large-$N$ confidence intervals based on a normal approximation.
We propose 
ratio-estimator-based and covariate-adjusted refinements for increased efficiency. We assess finite-sample empirical performance with a naturalistic simulation on real-world school social network data. The results demonstrate the reliability of the proposed methods in a realistic sample.

We then present our analysis of a field experiment on the effects of a program meant to promote anti-conflict norms and behavior among middle school students. In the experiment, schools were first randomly assigned to host the anti-conflict program and then sets of students within the host schools were randomly assigned to participate directly in the program.  The goal was to understand how attitudinal and behavioral effects on participants might transmit through their social network and affect their peers' behavior.

\section{Related Literature}

Our 
framework extends from
the foundational work of \citet{hudgens_halloran08}, who study two-stage, hierarchical randomized trials in which some groups are randomly assigned to host treatments; treatments are then assigned at random to units within the selected groups, and interference is presumed to operate only within groups. 
 Hudgens and Halloran provide randomization-based estimators for group-average causal effects, conditional on assignment strategies that determine the density of treatment within groups.  \citet{tchetgen_vanderweele2010} extend Hudgens and Halloran's results, providing conservative variance estimators, a framework for finite sample inference with binary outcomes, and extensions to observational studies. 
\citet{liu-hudgens2014-interference-intervals} develop asymptotic results for such two-stage designs.
Related to these contributions is work by \citet{rosenbaum07_interference}, which provides methods for inference with exact tests under partial interference. Under hierarchical treatment assignment and partial interference, estimation and inference can proceed assuming independence across groups.
In some settings, however, the hierarchical structuring may not be valid, as with experiments carried out over networks of actors that share links as a result of a complex, endogenous process.  
\citet{bowers-etal2013-interference} apply exact tests to evaluate parameters in models of spillover processes. Such a testing approach differs in its aims from ours, which focuses on estimating averages of potentially heterogenous unit-level causal effects.

A key contribution of this paper is to go beyond
the setting of hierarchical experiments with partial interference,
and to generalize estimation and inference theory to settings that exhibit arbitrary forms of interference and treatment assignment dependencies. In addition, our 
framework allows the analyst to work with different estimands, including both
types of group-average causal effects defined by the authors above  
as well as
average unit-level causal effects.  
Average unit-level causal effects are often the estimand of primary interest, as is the case, for example, when exploring unit-level characteristics that moderate the magnitude of treatment effects. 

\section{Treatment Assignment and Exposure Mappings\label{sec:exposure_mapping}}

In this section, we define the first two components of our analytical framework: the experiment design and exposure mapping.  We focus on the case of a randomized experiment with an arbitrary but known exposure mapping.  The first step is to distinguish between (i) treatment assignments over the set of experimental units and (ii) each unit's treatment exposure under a given assignment.  Treatment assignments can be manipulated arbitrarily with the experimental design.  However, treatment-induced exposures may be constrained on the basis of the varying potential for interference of different experimental units.  For example, interference or spillover effects may spread over a spatial gradient. If so, different treatment assignments may result in different patterns of interference depending on where treatments are applied on the spatial plane.

Formally, suppose we have a finite population $U$ of units indexed by $i=1,...,N$ on which a randomized experiment is performed.  Define a treatment assignment vector, $\z = (z_1, ..., z_N)'$, where $z_i \in \{1,..,M\}$ specifies which of $M$ possible treatment values that unit $i$ receives.  An {\it experimental design} contains a plan for randomly selecting a particular value of $\z$ from the $M^N$ different possibilities with predetermined probability $p_\z$.  Restricting our attention only to treatment assignments that can be generated by a given experimental design, define $\Omega = \{\z: p_\z > 0\}$, so that $\Z = (Z_1, ..., Z_N)'$ is a random vector with support $\Omega$ and $\Pr(\Z = \z) = p_\z$.  Our analysis below focuses on the case where the design is known in the following sense:
$\Pr(\Z = \z)$ for all $\z \in \Omega$ is known.  

We define 
an {\it exposure mapping} as 
a function that maps an assignment vector and unit-specific traits to an exposure value: $f : \Omega \times \Theta \rightarrow \Delta$, where $\theta_i \in \Theta$ quantifies relevant traits of unit $i$. The exposure mapping construction is functionally equivalent to the ``effective treatments'' function used by \citet{manski2012_identification_social}, 
though we find it helpful to denote separately the unit-specific attributes, $\theta_i$, that feed into the exposure mapping, $f(\cdot)$.
In applications we consider below, $\theta_i$ is unit $i$'s row in a network adjacency matrix.  More complex exposure mappings could take in $\theta_i$s that encode other traits of units and their peers---not only network ties, but also differences in age, gender, or other unit-level characteristics.  Or, $\theta_i$ could encode not only first-degree peer connections, but second-degree connections, third-degree, and so on.
 The codomain $\Delta$ contains all of the possible treatment-induced exposures 
that might be induced in the experiment.  The contents of $\Delta$ depend on the nature of interference.  These exposures may be represented as vectors, discrete classes, or 
scalar values.
As we will show formally below, each of the distinct exposures in $\Delta$ may give rise to distinct potential outcomes for each unit in $U$.  The estimation of causal effects under interference amounts to using information about {\it treatment assignments}, which come from the experiment's design, to estimate effects 
defined in terms
of {\it treatment-induced exposures}, which result from the interaction of the design (captured by $\Z$) and other underlying features of the population (captured by $f$ and the $\theta_i$s).  

To make things more concrete, consider some examples of exposure mappings.
The Neyman-Rubin causal model typically considers inference under an 
exposure mapping 
in which we set $\Delta = \{1,...,M\}$ and $f(\z, \theta_i) =f(\z) = z_i$ for all $i$.  This model has been a workhorse 
for much of the causal inference literature
\citep{neyman23, rubin78, holland86, imbens_rubin11}. An exposure mapping that allowed for completely arbitrary interference would be one for which $|\Delta| = |\Omega| \times N$, in which case each unit has a unique type of exposure under each treatment assignment, and $f(\z, \theta_i)$ would be unique for each $\z$. 
If such an exposure mapping were valid, then it is clear that there would be no meaningful way to use the results of the experiment to estimate average exposure-specific effects (although other types of causal effects may admit well-behaved estimators). Instead, the analyst must use substantive judgment 
about the extent of interference
to fix a 
mapping 
somewhere between the traditional randomized experiment and completely arbitrary exposure mappings in order to carry out analyses under interference. For example, \citet{hudgens_halloran08} consider a setting that allows unit $i$'s exposure to vary with each possible treatment assignment within $i$'s group, but, where conditional on the assignment for $i$'s group, $i$'s exposure does not vary in the treatment assignments of other groups.  Then, $\theta_i$ would be unit $i$'s group index, $|\Delta|$ would equal the largest number of assignment possibilities for any group, and $f(.)$ would map each assignment possibility for units in unit $i$'s group to a separate exposure condition.  The types of effects that \citet{hudgens_halloran08} construct are ones that average over these exposures for each unit.  
Below we discuss implications of using an exposure model that does not fully account for interference, drawing
connections back to the estimators in \citet{hudgens_halloran08}.  In the simulation study below and in the application, we provide more examples of exposure mappings.  Finally, our characterization of exposures is ``reduced form'' in that it does not distinguish between the channels through which interference occurs \citep{ogburn-venderweele2014}. The exposure mapping does not distinguish between effects that emanate directly from treatments being assigned to peers or that are mediated by changes in peers' outcomes (\citealp{eckles-etal2014}, pp. 8-9;  \citealp{manski95_identification_book}, Ch. 7).

Units' probabilities of falling into one or another exposure condition are crucial for the estimation strategy that we develop below. Define the exposure that unit $i$ receives as $D_i = f(\Z, \theta_i)$, a random variable with support $\Delta_i \subseteq \Delta$ and for which $\Pr(D_i = d) = \pi_i(d)$. Note that because $|\Delta| \le |\Omega| \times N$, $\Delta$ is a finite set of $K \le |\Omega|\times N$ values, such that $\Delta=\{d_1, ..., d_K\}$.  Then for each unit, $i$, we have a vector of probabilities, $(\pi_i(d_1),...,\pi_i(d_K))' = \bpi_i$.  Invoking \citet{imbens00}'s {\it generalized propensity score}, we call $\bpi_i$ the {\it generalized probability of exposure} for $i$. A unit $i$'s generalized probability of exposure tells us the probability of $i$ being subject to each of the possible exposures in $\{d_1,...,d_K\}$. We have 
$$
\pi_i(d_k) = \sum_{\z \in \Omega} \I(f(\z,\theta_i)=d_k)\Pr(\Z=\z) = {\sum_{\z \in \Omega} p_\z \I(f(\z,\theta_i)=d_k)}.
$$
Thus the generalized probability of exposure for unit $i$ is
also
known exactly. Each component probability, $\pi_i(d_k)$, is equal to the expected proportion of treatment assignments that induce exposure $d_k$ for unit $i$.

Below, we will refer to joint exposure probabilities when discussing variance estimators.  That is, we define $\pi_{ij}(d_k)$ as the probability of the joint event that both units $i$ and $j$ are subject to exposure $d_k$, and we define $\pi_{ij}(d_k, d_l)$ as the probability of the joint event that units $i$ and $j$ are subject to exposures $d_k$ and $d_l$, respectively.  To compute both individual and joint exposure probabilities from the experiment's design, first define the $N \times |\Omega|$ matrix
$$
\begin{array}{l}
\I_{k} = [\I(f(\z,\theta_i)=d_k)]_{\stackrel{\z \in \Omega}{i = 1,...,N}} = \\
\\
 \qquad \qquad \left[ \begin{array}{cccc}  \I(f(\z_1,\theta_1)=d_k)  &   \I(f(\z_2,\theta_1)=d_k)  & \hdots & \I(f(\z_{|\Omega|},\theta_1)=d_k) \\
  \I(f(\z_1,\theta_2)=d_k) &    \I(f(\z_2,\theta_2)=d_k) & \hdots & \I(f(\z_{|\Omega|},\theta_2)=d_k)  \\ \vdots & \vdots & \ddots & \\ \I(f(\z_1,\theta_N)=d_k)  &   \I(f(\z_2,\theta_N)=d_k)  & & \I(f(\z_{|\Omega|},\theta_N)=d_k) \end{array} \right],
  \end{array}
$$
which is a matrix of indicators for whether units are in exposure condition $k$ over possible assignment vectors.  Define the $|\Omega| \times |\Omega|$ diagonal matrix $\bP = \text{diag}(p_{\z_1}, p_{\z_2}, ..., p_{\z_{|\Omega|}})$. Then
$$
\I_k \bP \I_k' = \left[\begin{array}{cccc} \pi_1(d_k) & \pi_{12}(d_k) & \hdots & \pi_{1N}(d_k) \\ \pi_{21}(d_k) &  \pi_2(d_k) & \hdots & \pi_{2N}(d_k) \\ \vdots & \vdots & \ddots & \\  \pi_{N1}(d_k) & \pi_{N2}(d_k) &  &  \pi_N(d_k) \\ \end{array} \right],
$$
is an $N \times N$ symmetric matrix with individual exposure probabilities, the $\pi_{i}(d_k)$s, on the diagonal and joint exposure probabilities, the $\pi_{ij}(d_k)$s, on the off-diagonals.  The non-symmetric $N \times N$ matrix
$$
\I_k \bP \I_l' = \left[\begin{array}{cccc} 0 & \pi_{12}(d_k,d_l) & \hdots & \pi_{1N}(d_k,d_l) \\ \pi_{21}(d_k,d_l) & 0 & \hdots & \pi_{2N}(d_k,d_l) \\ \vdots & \vdots & \ddots & \\  \pi_{N1}(d_k,d_l) & \pi_{N2}(d_k,d_l) &  & 0 \\ \end{array} \right],
$$
yields all joint probabilities across exposure conditions $k$ and $l$.  The zeroes on the diagonal are due to the fact that a unit cannot be subject to multiple exposure conditions at once.

In practice, $|\Omega|$ may be so large that it is impractical to construct $\Omega$ to compute the $\bpi_i$s and the joint probability matrices exactly.  One may nonetheless approximate the $\bpi_i$s and joint probabilities with arbitrary precision through simulation; that is, produce $R$ random replicate $\z$s based on the randomization plan.  From these $R$ replicates, we can construct an $N \times R$ indicator matrix, $\widehat{\I}_k$, for each of the $k=1,...,K$ exposure conditions.  
Then, an estimator for $\mathbf{I}_k\mathbf{P}\mathbf{I}_k'$ that incorporates mild additive smoothing to ensure non-zero marginal probability estimates is $(\hat{\mathbf{I}}_k\hat{\mathbf{I}}_k' + \iota_N)/(R+1)$, where $\iota_N$ is an $N \times N$ identity matrix.  Similarly, an estimator for $\mathbf{I}_k\mathbf{P}\mathbf{I}_l'$, which does not admit additive smoothing due to zero joint inclusion probabilities, is $(\hat{\mathbf{I}}_k\hat{\mathbf{I}}_l')/R$.

\begin{prop} \label{prop:sim-converge}
As $R \rightarrow \infty$,
$$
(\hat{\mathbf{I}}_k\hat{\mathbf{I}}_k' + \iota_N)/(R+1) \overset{a.s.}{\rightarrow} \mathbf{I}_k\mathbf{P}\mathbf{I}_k', \text{ and } (\hat{\mathbf{I}}_k\hat{\mathbf{I}}_l')/R \overset{a.s.}{\rightarrow} \mathbf{I}_k\mathbf{P}\mathbf{I}_l'.
$$
\end{prop}
\noindent All proofs appear in the appendix.  Rates of convergence of $\widehat{\I}_k\widehat{\I}'_k/R$ are discussed in \citet{fattorini06} and \citet{aronowdiss}. 
Below we give guidance on selecting a value of $R$ based on a bound on the relative bias for an estimator of a target quantity.


\section{Average Potential Outcomes and Causal Effects}\label{avgpo}

We develop the case of estimating average unit-level causal effects of exposures.  An average unit-level causal effect is defined in terms of a difference between the average of units' potential outcomes under one exposure versus the average under another exposure.  
The starting point is the estimation of
average potential outcomes under each of the exposure conditions.
With that, the analyst is in principle free to compute a variety of causal quantities of interest, 
not just average unit-level causal effects.
For example, one could consider effects that are defined as differences between the average of potential outcomes under one set of exposures versus the average under another set of exposures. 
The 
direct, indirect and overall effects of \citet{hudgens_halloran08}
are defined in this way using the construction of the ``individual average potential outcome.'' The hierarchical designs that they consider are specifically tailored to ensure that estimators for such effects are non-parametrically identified. 
While our focus is on estimating exposure-specific causal effects that are 
defined for arbitrary designs,
such design-specific estimators can 
certainly be derived and analyzed using the framework developed here.
 Our focus on the average of unit-level, exposure-specific causal effects is due to its being the natural extension of the ``average treatment effect'' that is the focus of much current causal inference and program evaluation literature (e.g., \citealp{imbens_wooldridge09}).  

Suppose 
all units have non-zero probabilities of being subject to each of the $K$ exposures:
$0 < \pi_i(d_k) < 1$ for all $i$ and $k$. (When $\pi_i(d_k) = 0$ for some units, then design-based estimation of average potential outcomes and causal effects must be restricted to the subset of units for which $\pi_i(d_k) > 0$.) 
In the most general terms, each $\z \in \Omega$ can generate a potential outcome for unit $i$. We label the randomization potential outcome of unit $i$ associated with $\z$ as  $y_i^r(\z)$.  These randomization potential outcomes are fixed for all units in the population and do not depend on the value of randomized treatment, $\Z$.
A condition that we use in our analysis below is that the exposure mapping fully characterizes interference:
\begin{cond}[Properly specified exposure mapping]
For all $i \in \{1,...,N\}$ and $\z, \z' \in \Omega$ such that $f(\z,\theta_i) = f(\z',\theta_i)$, $y_i^r(\z) = y_i^r(\z')$. 
\label{cond:sutva}
\end{cond} 
\noindent Given Condition \ref{cond:sutva}, each unit $i$ has $|\Delta| = K$ potential outcomes, which we can write in terms of the exposure conditions as $(y_i(d_1),...,y_i(d_K))$, where $y_i(d_k) = y_i^r(\z), \forall i \in \{1,...,N\}$, $k \in \{1,...,K\}$, and $\z \in \Omega$ such that $f(\z,\theta_i) = d_k$.

Let $Y_i$ be the observed outcome for unit $i$. We assume the following consistency condition that relates the observed data to potential outcomes under the exposure model \citep{vanderweele2009-consistency}:
\begin{cond}[Consistent potential outcomes]
$$
Y_i = \sum_{k=1}^K \mathbf{I}(D_i = d_k) y_i(d_k), \forall i \in \{ 1,...,N \}. $$
\label{cond:consistency}
\end{cond} 
\noindent 
\noindent Although 
SUTVA
is violated at the level of treatment assignment (i.e., the individual $z_i$ values), 
Conditions \ref{cond:sutva} and \ref{cond:consistency}
restore a form of SUTVA with respect to exposures, in a manner that is conceptually similar to \citet{hudgens_halloran08}'s stratified interference assumption. Below we examine implications of violating Conditions \ref{cond:sutva} and \ref{cond:consistency}---i.e., implications of misspecifying the exposure mapping. Throughout, unless otherwise noted, we will be assuming Conditions \ref{cond:sutva} and \ref{cond:consistency}.

We seek estimates for all $k$ of $\mu(d_k) = \frac{1}{N}\sum_{i=1}^N y_i(d_k)=\frac{1}{N}y^T(d_k)$, where $y^T(d_k)$ is the total of the potential outcomes under $d_k$.  This allows us to define an average causal effect of being in exposure condition $d_k$ as compared to being in condition $d_l$ as
$$ 
\tau(d_k, d_l) = \frac{1}{N}\sum_{i=1}^N y_i(d_k) - \frac{1}{N}\sum_{i=1}^N y_i(d_l).
$$

The number of units in the population, $N$, is fixed, but we cannot estimate $y^T(d_k)$ directly, as we only observe $y_i(d_k)$ for those with $D_i = d_k$.  However, by design, the collection of units for which we observe $y_i(d_k)$ is an unequal-probability without-replacement sample from $(y_1(d_k),...,y_N(d_k))$, with the sampling probabilities known exactly.  By \citet{horvitz_thompson}, a design-based estimator for $y^T(d_k)$ is the inverse probability weighted estimator
\begin{equation}
\widehat{y^T_{HT}}(d_k) = \sum_{i=1}^N \I(D_i=d_k)\frac{Y_i}{\pi_i(d_k)}
\label{eq:ht_esimator}
\end{equation}
Below, we consider variance-reducing refinements to this estimator. We start with an analysis of $\widehat{y^T_{HT}}(d_k)$ because it very clearly reveals first-order issues for estimation under interference.
Estimator \ref{eq:ht_esimator} is unbiased, and its variance is characterized in Lemma \ref{lemma:ht}.
\begin{lemma}\label{lemma:ht}
$$\E[\widehat{y^T_{HT}}(d_k)] = \sum_{i=1}^N y_i(d_k)$$
\begin{align} 
\Var[\widehat{y^T_{HT}}(d_k)]  
& = \sum_{i=1}^N \pi_i(d_k)[1-\pi_i(d_k)] \left[ \frac{y_i(d_k)}{\pi_i(d_k)} \right]^2 \nonumber \\ & \hspace{1em} + \sum_{i=1}^N \sum_{j \ne i}  [\pi_{ij}(d_k)- \pi_i(d_k)\pi_j(d_k)]\frac{y_i(d_k)}{\pi_i(d_k)} \frac{y_j(d_k)}{\pi_j(d_k)}. \label{eq:total_variance}
\end{align}
\end{lemma}
Above we indicated that one can approximate $\mathbf{I}_k\mathbf{P}\mathbf{I}_k'$ with $(\hat{\mathbf{I}}_k\hat{\mathbf{I}}_k' + \iota_N)/(R+1)$, which has diagonal elements,
$$
\widehat{\pi_i}(d_k) = \frac{X_i + 1}{R + 1},
$$
where $X_i = \sum_{r=1}^R \mathbf{I}_r(f(\mathbf{z}_r, \theta_i) = d_k)$ and $r = 1,...,R$ indexes the replicates. Define the Horvitz-Thompson estimator that uses the $\widehat{\pi_i}(d_k)$ estimates: 
$$
\widehat{y^T_{HT,R}} (d_k) = \sum_{i=1}^N \mathbf{I}(D_i = d_k)\frac{Y_i}{\widehat{\pi_i}(d_k)}.
$$
Following  \citet{fattorini06}, the following proposition provides guidance on choosing $R$ in terms of a bound on the relative bias for estimating $\widehat{y^T} (d_k)$.

\begin{prop} \label{prop:bias-bound}
The relative bias for $\widehat{y^T_{HT,R}} (d_k)$ is bounded as  
$$
\left| \frac{\E[\widehat{y^T_{HT,R}} (d_k)] - y^T(d_k)}{y^T(d_k)} \right| \le (1-\pi_0(d_k))^{R+1},
$$
where $\pi_0(d_k) = \min_i\{ \pi_i(d_k) \}$.
\end{prop}
For a relative bias target of $b$ and given some approximation of $\pi_0(d_k)$, the bound implies selecting a number of replicates $
R \ge \log(b)/\log(1-\pi_0(d_k)) - 1$.  Thus, for $b = .005$ and $\pi_0(d_k) = .0005$, this would imply at least 10,593 replicates.   
Given the apparent computational feasibility of producing enough replicates so as to render relative biases negligible, from here on our analysis  assumes that we are working with $\mathbf{I}_k\mathbf{P}\mathbf{I}_k'$ and $\mathbf{I}_k\mathbf{P}\mathbf{I}_l'$.

Given the estimator of the total of the $N$ potential outcomes under exposure $d_k$, a natural estimator for the mean is thus $
\widehat{\mu_{HT}}(d_k) = (1/N)\widehat{y^T_{HT}}(d_k), \label{eq:exposure_mean}
$
with variance 
$
\Var(\widehat{\mu_{HT}}(d_k)) = (1/N^2)\Var[\widehat{y^T_{HT}}(d_k)].\label{eq:var_of_exposure_mean}
$ 
This allows us to construct the difference in estimated means
\begin{equation}
\widehat {\tau_{HT}}(d_k,d_l) = \widehat{\mu_{HT}}(d_k) - \widehat{\mu_{HT}}(d_l) = \frac{1}{N}\left[\widehat{y^T_{HT}}(d_k) - \widehat{y^T_{HT}}(d_l)\right] \label{eq:ht_causal_effect}
\end{equation}
which is an estimator of $\tau(d_k,d_l) = \frac{1}{N}\sum_{i=1}^N\left[y_i(d_k)-y_i(d_l)\right]$, the average unit-level causal effect of 
exposure $k$ versus 
exposure $l$.  

\begin{prop} \label{prop:ht}
\begin{equation}
\E[\widehat {\tau_{HT}}(d_k,d_l)] = \frac{1}{N}\sum_{i=1}^N y_i(d_k) - \frac{1}{N}\sum_{i=1}^N y_i(d_l)
\label{eq:etauhat}
\end{equation}
\begin{align}
\Var(\widehat {\tau_{HT}}(d_k,d_l)) = & \frac{1}{N^2}\left\{\Var[\widehat{y^T_{HT}}(d_k)] + \Var[\widehat{y^T_{HT}}(d_l)] \right. \nonumber \\ & \left. \hspace{3em} - 2\Cov[\widehat{y^T_{HT}}(d_k),\widehat{y^T_{HT}}(d_l)]\right\},\label{eq:tru_var}
\end{align}
where 
\begin{align}
\Cov[\widehat{y^T_{HT}}(d_k),\widehat{y^T_{HT}}(d_l)]
& = \sum_{i=1}^N \sum_{j \ne i}  \frac{y_i(d_k)}{\pi_i(d_k)} \frac{y_j(d_l)}{\pi_j(d_l)} \left[\pi_{ij}(d_k,d_l)- \pi_i(d_k)\pi_j(d_l) \right] \nonumber \\ & \hspace{1em} - \sum_{i=1}^N y_i(d_k)y_i(d_l). \label{eq:totals_covariance}
\end{align}
\end{prop}
Expressions \eqref{eq:total_variance} and \eqref{eq:totals_covariance} allow us to see the conditions under which exact variances are identified.  So long as all joint exposure probabilities are non-zero (that is, $\pi_{ij}(d_k) > 0$ for all $i,j$), unbiased estimators for $\Var[\widehat{y^T_{HT}}(d_k)]$ are identified for the population $U$.  Because we only observe one potential outcome for each unit, the last sum in \eqref{eq:totals_covariance} is always unidentified, and thus $\Cov[\widehat{y^T_{HT}}(d_k),\widehat{y^T_{HT}}(d_l)]$ is always unidentified.  This is a familiar problem in estimating the randomization variance for the average treatment effect---e.g., \citet{neyman23} or \citet[A32-A34]{freedman_pisani_purves98}.  If $\pi_{ij}(d_k) = 0$ for some $i,j$, $\Var[\widehat{y^T_{HT}}(d_k)]$ is unidentified. Similarly, if $\pi_{ij}(d_k,d_l) = 0$ for some $i,j$, then additional components of $\Cov[\widehat{y^T_{HT}}(d_k),\widehat{y^T_{HT}}(d_l)]$ are unidentified. 
Nonetheless, we can always identify estimators for $\Var[\widehat{y^T_{HT}}(d_k)]$ and $\Cov[\widehat{y^T_{HT}}(d_k),\widehat{y^T_{HT}}(d_l)]$ that are guaranteed to have nonnegative bias.  Thus, we can always identify a conservative approximation to the exact variances.  We discuss this in the next section. 

\section{Variance Estimators}
\label{varests}

We derive conservative estimators for both $\Var[\widehat{y^T_{HT}}(d_k)]$ and $\Var[\widehat {\tau_{HT}}(d_k,d_l)]$. The formulations in this section follow from \citet{aronow_samii_zeropairwise} which considers conservative variance estimation for generic sampling designs with some zero pairwise inclusion probabilities.  Although not necessarily unbiased, the estimators we present here are guaranteed to have a nonnegative bias relative to the randomization distributions of the estimators.

Given $\pi_{ij}(d_k) > 0$ for all $i,j$, the Horvitz-Thompson estimator for $\Var[\widehat{y^T_{HT}}(d_k)]$ is
\begin{align}
\widehat{\Var}[\widehat{y^T_{HT}}(d_k)] 
& = \sum_{i \in U}\I(D_i=d_k)[1-\pi_{i}(d_k)] \left[ \frac{Y_i}{\pi_i(d_k)} \right]^2 \nonumber \\ & \hspace{1em} + \sum_{i \in U} \sum_{j \in U \backslash i} \I(D_i=d_k)\I(D_j=d_k) \nonumber \\ & \hspace{4em} \times \frac{\pi_{ij}(d_k)-\pi_{i}(d_k)\pi_{j}(d_k)}{\pi_{ij}(d_k)}\frac{Y_i}{\pi_i(d_k)} \frac{Y_j}{\pi_j(d_k)}. \label{eq:ht_variance_estimator}
\end{align}

\begin{lemma}\label{lem:varun}
If $\pi_{ij}(d_k) > 0$  for all $i,j$, then 
$
\E[\widehat{\Var}[\widehat{y^T_{HT}}(d_k)] ] = \Var[\widehat{y^T_{HT}}(d_k)].
$
\end{lemma}
\noindent Lemma \ref{lem:varun} follows from unbiasedness of the Horvitz-Thomspon estimator for measurable designs. Then an unbiased estimator for the variance of $\widehat{\mu_{HT}}(d_k)$ is 
$
\widehat{\Var}[\widehat{\mu_{HT}}(d_k)] = (1/N^2)\widehat{\Var}[\widehat{y^T_{HT}}(d_k)].
$

In the case where $\pi_{ij}(d_k) = 0$ for some $i,j$, the Horvitz-Thompson estimator of $\Var[\widehat{y^T_{HT}}(d_k)]$ is not unbiased, but its bias is readily characterized.  
\begin{prop}\label{lem:varbias}
If $\pi_{ij}(d_k) = 0$  for some $i,j$, then 
$
\E[\widehat{\Var}[\widehat{y^T_{HT}}(d_k)] ] = \Var[\widehat{y^T_{HT}}(d_k)] + A,
$
where $$
A= \sum_{i \in U} \sum_{j \in \{U \backslash i:\pi_{ij}(d_k)=0\}} y_i(d_k)y_j(d_k).
$$
\end{prop}
\noindent A proof for Proposition \ref{lem:varbias} follows from \citet[Proposition 1, reproduced in the appendix below]{aronow_samii_zeropairwise}.

Note that $\widehat{\Var}[\widehat{\mu_{HT}}(d_k)]$ is guaranteed to have nonnegative bias when $y_i(d_k)y_j(d_k) \ge 0$ for all $i,j$ with $\pi_{ij}(d_k) = 0$.  The bias will be small when the terms in the sum tend to offset each other, as when the relevant $y_i(d_k)$ and $y_j(d_k)$ values are centered on 0 and have low correlation with each other. (This notation requires that we define ${0}/{0} = 0$.) 

 Another option is to use the following correction term (derived via Young's inequality),
$$
\widehat{A_2}(d_k) = \sum_{i \in U}\sum_{j \in \{ U \backslash i : \pi_{ij}(d_k) = 0\}} \left[\frac{\I(D_i=d_k)Y_i^2}{2\pi_{i}(d_k)} + \frac{\I(D_j=d_k)Y_j^2}{2\pi_{j}(d_k)}\right],
$$
noting that $\widehat{A_2}(d_k) = 0$ if $\pi_{ij}(d_k) > 0 \textrm{ for all } i,j$. 
\begin{prop}\label{a2}
$$
\E \left[ \widehat{\Var}[\widehat{y^T_{HT}}(d_k)] + \widehat{A_2}(d_k) \right] \ge \Var[\widehat{y^T_{HT}}(d_k)],
$$
\end{prop}
\noindent A proof for Propositon \ref{a2} follows directly from \citet[Corollary 2, reproduced in the appendix below]{aronow_samii_zeropairwise}. Then let
$
\widehat{\Var_A}[\widehat{\mu_{HT}}(d_k)] = (1/N^2)\left[\widehat{\Var}[\widehat{y^T_{HT}}(d_k)] + \widehat{A_2}(d_k)\right].
$ $\widehat{\Var_A}[\widehat{\mu_{HT}}(d_k)]$ then provides a conservative estimator for the variance of the estimated average of potential outcomes under exposure $d_k$.  

As discussed above, $\Cov[\widehat{y^T_{HT}}(d_k),\widehat{y^T_{HT}}(d_l)]$ is unidentified, which is to say that there exist no unbiased or consistent estimators for this quantity.  However, we can compute an approximation that is guaranteed to have expectation less than or equal to the true covariance, providing a conservative (here, nonnegatively biased) estimator for $\Var(\widehat {\tau_{HT}}(d_k,d_l))$.  For the case where $\pi_{ij}(d_k,d_l) > 0$ for all $i,j$ such that $i \neq j$, we propose the Horvitz-Thompson-type estimator for the covariance
\begin{align}
\widehat{\Cov}[\widehat{y^T_{HT}}(d_k),\widehat{y^T_{HT}}(d_l)]  & = \sum_{i \in U} \sum_{j \in U \backslash i}\Bigg[ \frac{\I(D_i = d_k)\I(D_j = d_l)}{\pi_{ij}(d_k,d_l)} \frac{Y_i}{\pi_i(d_k)} \frac{Y_j}{\pi_j(d_l)} \nonumber \\
& \hspace{5.5em} \times [\pi_{ij}(d_k,d_l) - \pi_i(d_k)\pi_j(d_l)]\Bigg] 
\nonumber \\ & \hspace{1em} - \sum_{i \in U}\left[ \frac{\I(D_i = d_k) Y_i^2}{2\pi_i(d_k)} + \frac{\I(D_i = d_l) Y_i^2}{2\pi_i(d_l)}\right]. \label{eq:ht_cov_estimator}
\end{align}

\begin{prop}\label{ncov}
If $\pi_{ij}(d_k,d_l) > 0$ for all $i,j$ such that $i \neq j$,
$$
\E \left[ \widehat{\Cov}[\widehat{y^T_{HT}}(d_k),\widehat{y^T_{HT}}(d_l)]  \right] \le {\Cov}[\widehat{y^T_{HT}}(d_k),\widehat{y^T_{HT}}(d_l)],
$$
\end{prop}
\noindent A proof for Proposition \ref{ncov} follows from noting that the term on the second line in expression \eqref{eq:ht_cov_estimator} has expected value less than or equal to the quantity in the last line of expression \eqref{eq:totals_covariance}, again via Young's inequality. See \citet[Proposition 2, reproduced in the appendix below]{aronow_samii_zeropairwise} for greater detail.

$\widehat{\Cov}[\widehat{y^T_{HT}}(d_k),\widehat{y^T_{HT}}(d_l)]$ is exactly unbiased if,  for all $i \in U$, $y_i(d_l) = y_i(d_k)$, implying no effect associated with condition $l$ relative to condition $k$. 
\begin{prop}\label{prop:no-bias-cov}
If $\pi_{ij}(d_k,d_l) > 0$ for all $i,j$ such that $i \neq j$ and for all $i \in U$, $y_i(d_l) = y_i(d_k)$
$$
\E \left[ \widehat{\Cov}[\widehat{y^T_{HT}}(d_k),\widehat{y^T_{HT}}(d_l)]  \right] = {\Cov}[\widehat{y^T_{HT}}(d_k),\widehat{y^T_{HT}}(d_l)],
$$
\end{prop}
\noindent A proof follows from \citet[Corollary 1, reproduced in the appendix below]{aronow_samii_zeropairwise}.

For the case where $\pi_{ij}(d_k,d_l) = 0$ for some $i,j$ and $k,l$, we can refine the expression for the covariance given in \eqref{eq:totals_covariance} to 
\begin{align}
\Cov[\widehat{y^T_{HT}}(d_k),\widehat{y^T_{HT}}(d_l)] = & \sum_{i \in U} \sum_{j \in \{U \backslash i : \pi_{ij}(d_k,d_l) > 0\}} \Bigg[ \frac{y_i(d_k)}{\pi_i(d_k)}\frac{y_j(d_l)}{\pi_j(d_l)}\nonumber \\ & \hspace{4em} \times [\pi_{ij}(d_k,d_l) - \pi_i(d_k)\pi_j(d_l)] \Bigg] \nonumber \\ & \hspace{1em} - \sum_{i \in U} \sum_{j \in \{U : \pi_{ij}(d_k,d_l) = 0\}}y_i(d_k) y_j(d_l),\label{eq:totals_covariance_general}
\end{align}
where the term on the last line subsumes the term on the last line in expression \eqref{eq:totals_covariance}.  This leads us to propose a more general estimator for the covariance
\begin{align}
\widehat{\Cov_A}[\widehat{y^T_{HT}}(d_k),\widehat{y^T_{HT}}(d_l)]  & = \sum_{i \in U} \sum_{j \in \{ U \backslash i : \pi_{ij}(d_k,d_l) > 0 \}}\Bigg[ \frac{\I(D_i = d_k)\I(D_j = d_l)}{\pi_{ij}(d_k,d_l)} \nonumber \\ & \hspace{10em} \times \frac{Y_i}{\pi_i(d_k)} \frac{Y_j}{\pi_j(d_l)}  \nonumber \\
& \hspace{10em} \times [\pi_{ij}(d_k,d_l) - \pi_i(d_k)\pi_j(d_l)] \Bigg] \nonumber  \\  & \hspace{1em} - \sum_{i \in U} \sum_{j \in \{ U : \pi_{ij}(d_k,d_l) = 0\} }\left[ \frac{\I(D_i = d_k) Y_i^2}{2\pi_i(d_k)} \right. \nonumber \\ & \hspace{10em} \left. + \frac{\I(D_j = d_l) Y_j^2}{2\pi_j(d_l)}\right]. \label{eq:ht_cov_general_estimator}
\end{align}
\begin{prop}\label{cova}
$$
\E \left[ \widehat{\Cov_A}[\widehat{y^T_{HT}}(d_k),\widehat{y^T_{HT}}(d_l)]  \right] \leq {\Cov}[\widehat{y^T_{HT}}(d_k),\widehat{y^T_{HT}}(d_l)],
$$
\end{prop}
\noindent A proof again follows from the fact the term in the last line in \eqref{eq:ht_cov_general_estimator} has expected value no greater than the term in the last line of \eqref{eq:totals_covariance_general} by Young's inequality.

Based on the variance expressions and correction terms defined above, we obtain a conservative variance estimator for $\Var(\widehat {\tau_{HT}}(d_k,d_l))$ as 
\begin{align}
\widehat{\Var}[\widehat {\tau_{HT}}(d_k,d_l)] & = \frac{1}{N^2} \left\{ \widehat{\Var}[\widehat{y^T_{HT}}(d_k)] + \widehat{A_2}(d_k) + \widehat{\Var}[\widehat{y^T_{HT}}(d_l)] + \widehat{A_2}(d_l) \right.\nonumber \\ & \left.\hspace{3.75em} - 2\widehat{\Cov_A}[\widehat{y^T_{HT}}(d_k),\widehat{y^T_{HT}}(d_l)]
\right\} .\label{eq:ate_var_estimator}
\end{align}

\begin{prop}\label{consvar}
$$
\E \left[ \widehat{\Var}[\widehat {\tau_{HT}}(d_k,d_l)]   \right] \ge {\Var}[\widehat {\tau_{HT}}(d_k,d_l)] ,
$$
\end{prop}
\noindent The result follows  from Proposition \ref{a2}, Proposition \ref{cova}, and linearity of expectations.


\section{Asymptotics and Intervals}
\label{asymptot}


Consider a sequence of nested populations indexed by size $N$, $(U_N)$. To define a notion of asymptotic growth, we let $N$ tend to infinity, allowing for the experimental design to be reapplied anew to each $U_N$, subject to the conditions defined below  \citep{brewer1979,isaki_fuller82}. 
%
Consistency and the asymptotic validity of Wald-type confidence intervals will then follow from restrictions on the growth process of the design and exposure mapping.

\subsection{Consistency} \label{consist}

We first establish conditions for the estimator $\widehat{\tau_{HT}}(d_k,d_l)$ to converge to $\tau(d_k,d_l)$ as $N$ grows. We will show that, under two regularity conditions, $\widehat{\tau_{HT}}(d_k,d_l) - \tau(d_k,d_l) \overset{p}{\longrightarrow} 0$ as $N \rightarrow \infty$.

\begin{cond}[Boundedness of potential outcomes and exposure probabilities] 
Potential outcomes and exposure probabilities are bounded, so that, for all values $i$ and $d_k$, $|y_i(d_k)| \leq c_1  < \infty$ and $|1/\pi_i(d_k)| \leq c_2 < \infty$. \label{cond1}
\end{cond}
\noindent  Condition \ref{cond1} can be relaxed, though Condition \ref{cond2} would likely need to be strengthened accordingly.

We will also make an assumption about the amount of dependence in exposure conditions in the population. Define a pairwise dependency indicator $g_{ij} $ such that if $g_{ij} = 0$, then $D_i \independent D_j$, else let $g_{ij} = 1$. 

\begin{cond}[Restrictions on pairwise dependence] $\sum_{i=1}^{N} \sum_{j=1}^{N} g_{ij} = o(N^{2})$. \label{cond2}
\end{cond}
\noindent Condition \ref{cond2} entails that, as $N$ grows, the amount of pairwise clustering in exposure conditions induced by the design and exposure mapping is limited in scope.   As units are added to the sample, the number of new non-zero entries in the expanding pairwise correlation matrix of exposures should be limited by the order condition.


\begin{prop} \label{prop:tau-consistent}
Given Conditions \ref{cond1} and \ref{cond2}, 
$\widehat{\tau_{HT}}(d_k,d_l) - \tau(d_k,d_l)  \overset{p}{\longrightarrow} 0$ as $N \rightarrow \infty$. 
\end{prop}

\subsection{Confidence Intervals}

We now establish conditions for the asymptotic validity of Wald-type confidence intervals under stricter conditions on the asymptotic growth process. Consistency for the variance estimators, asymptotic normality, and therefore asymptotic validity of confidence intervals, follow straightforwardly when the amount of dependence across units in the population is limited. 

We shall assume that Condition \ref{cond1} holds, but will strengthen Condition \ref{cond2} to ensure that dependence across exposures is limited in scope. Unlike Condition \ref{cond2}, we will exploit joint independence of observations rather than pairwise independence. Define a binary dependency indicator $h_{ij}$ over all pairs $(i,j) \in \{1,2,...,N\} \times  \{1,2,...,N\}$, where $h_{ij}$ satisfies the following: for any pair of disjoint sets $\Gamma_1$ and $\Gamma_2 \subseteq \{1,...,N\}$ such that there exists no pair $(i,j)$ with $h_{ij} = 1$ and either (i) $i \in \Gamma_1$ and $j \not\in \Gamma_2$, (ii) $j \in \Gamma_1$ and $i \not\in \Gamma_2$, (iii) $i \not\in \Gamma_1$ and $j \in \Gamma_2$, or  (iv) $j \not\in \Gamma_1$ and $i \in \Gamma_2$, $\{D_i, i \in \Gamma_1\}$ and $\{D_i, i \in \Gamma_2\}$ are independent.

%
%
%
%
%
%
%


\begin{cond}[Local dependence] \label{cond:local-dep} There exists a finite constant $m$ such that, for all $U_N$, $i \in 1,...,N$, $\sum_{j=1}^{N} h_{ij} \leq m$. 
%
%
%
%
%
%
%
\end{cond}


\noindent Condition \ref{cond:local-dep} is equivalent to assuming that dependencies across exposures can be represented by a dependency graph such that the maximal degree of each unit tends to be limited relative to $N$. Condition \ref{cond:local-dep} will allow us to straightforwardly invoke a central limit theorem for random fields as derived via Stein's Method \citep[Example 2.4.1]{chenshao2004}. (The authors thank Betsy Ogburn for the suggestion of the use of Stein's Method in this setting.) Finite $m$ ensures that our variance estimators will converge at a sufficiently fast rate. Note that Condition \ref{cond:local-dep} subsumes Condition \ref{cond2}, as $\sum_{i=1}^{N} \sum_{j=1}^{N} g_{ij} = O(N)$ when Condition \ref{cond:local-dep} holds.

It is illustrative to consider settings where Condition \ref{cond:local-dep} holds. For Bernoulli-randomized designs, Condition \ref{cond:local-dep} would hold if interference were characterized by first-order dependence on a graph connecting units and network degrees were bounded above by some value $m$.  Condition \ref{cond:local-dep} also generalizes the partial interference setting considered by, e.g., \citet{sobel06} and \citet{liu-hudgens2014-interference-intervals} given finite subpopulations across which interference is localized (in this case, $m$ would be the size of the largest subpopulation). However, Condition \ref{cond:local-dep} would be violated if changing the treatment  assigned to one unit would affect the exposure received by all $N$ units.  In comparing Conditions \ref{cond2} and \ref{cond:local-dep}, note that Condition \ref{cond2} is a restriction on the order of growth of pairwise dependencies, while Condition \ref{cond:local-dep} requires local dependence.  The latter condition is more restrictive as it imposes conditions on all higher-order joint inclusion probabilities. It is possible that Condition \ref{cond2} could hold, but Condition \ref{cond:local-dep} would be violated, if, for example, there exists a single unit for which the number of associated pairwise dependencies tended to infinity in $N$.Ó

\begin{cond}[Nonzero limiting variance.] \label{cond:nonzero-var}
$N \Var [ \widehat {\tau_{HT}}(d_k,d_l)] \rightarrow c$, where $c > 0$.
\end{cond}

\noindent Convergence of $N \Var [ \widehat {\tau_{HT}}(d_k,d_l)]$ to a nonnegative constant is generally ensured by Conditions \ref{cond1} and \ref{cond:local-dep}, sufficient for root-$n$ consistency of $\widehat {\tau_{HT}}(d_k,d_l)$. Condition \ref{cond:nonzero-var} is a mild regularity condition that ensures that this constant is positive, and rules out degenerate cases (e.g., all outcomes are zero). 

\begin{prop}\label{prop:tau-interval}
Given Conditions \ref{cond1}, \ref{cond:local-dep},  and \ref{cond:nonzero-var},  Wald-type intervals constructed as
$$
\widehat {\tau_{HT}}(d_k,d_l) \pm z_{1-\alpha/2}\sqrt{\widehat{\Var}[\widehat {\tau_{HT}}(d_k,d_l)]}
$$
will tend to cover $\tau_{HT}(d_k,d_l)$ at least $100(1-\alpha)\%$ of the time for large $N$.
\end{prop}

\section{Refinements} \label{interfererefine}

The mean and difference-in-means estimators presented thus far are unbiased by sample theoretic arguments, and we have derived conservative variance estimators.  However, we may wish to improve efficiency by incorporating auxiliary covariate information.  In addition, by analogy to results from the unequal probability sampling literature, 
ratio approximations of the Horvitz-Thompson estimator
may significantly reduce mean squared error with little cost in terms of bias \citep[pp. 181-184]{sarndal_etal92}. 
We discuss such refinements here. 

\subsection{Covariance Adjustment}

Auxiliary covariate information may help to improve efficiency. A first method of covariance adjustment is based on the so-called ``difference estimator'' \citep[Ch. 6]{raj65, sarndal_etal92}.  Covariance adjustment of this variety can reduce the randomization variance of the estimated exposure means and average causal effects without compromising unbiasedness.  In addition, the difference estimator addresses the problem of location non-invariance that afflicts Horvitz-Thompson-type estimators \citep[pp. 9-10]{fuller09_samp_stat}. The estimator requires prior knowledge of how outcomes relate to covariates, perhaps obtained from analysis of auxiliary datasets. 

Assume an auxiliary covariate vector $\mathbf{x_i}$ is observed for each $i$.  We have some predefined function $g\left(\mathbf{x_i}, \mathbf{\xi_i}(d_k) \right) \rightarrow \mathbb{R}$, where $\mathbf{\xi_i}$ is a parameter vector. Ideally $g(\cdot)$ is calibrated on auxiliary data to produce values that approximate $y_i(d_k)$.  We assume $\Cov[g\left(\mathbf{x_i}, \mathbf{\xi_i}(d_k) \right),Z_i] = 0$ as a sufficient condition for unbiasedness. 
Define
\begin{equation}
\label{eq:gen1}
\widehat{y^{T}_{G}}(d_k) = \sum_{i=1}^N \I(D_i=d_k)\frac{Y_i}{\pi_i(d_k)}  - \sum_{i=1}^N \I(D_i=d_k)\frac{g\left(\mathbf{x_i}, \mathbf{\xi_i}(d_k) \right)}{\pi_i(d_k)} +  \sum\limits_{i=1}^Ng\left(\mathbf{x_i}, \mathbf{\xi_i}(d_k) \right),
\end{equation}
which is unbiased for $y^T(d_k)$ by
$$
\E\left[- \sum_{i=1}^N \I(D_i=d_k)\frac{g\left(\mathbf{x_i}, \mathbf{\xi_i}(d_k) \right)}{\pi_i(d_k)} +  \sum\limits_{i=1}^Ng\left(\mathbf{x_i}, \mathbf{\xi_i}(d_k) \right)\right] = 0.
$$
Define $\epsilon_i(d_k) = Y_i - g\left(\mathbf{x_i}, \mathbf{\xi_i}(d_k) \right)$ for cases with $D_i = d_k$. Then, by substitution,
\begin{equation}
\label{eq:gen2}
\widehat{y^{T}_{G}}(d_k) = \sum_{i=1}^N \I(D_i=d_k)\frac{\epsilon_i(d_k)}{\pi_i(d_k)}  +  \sum\limits_{i=1}^Ng\left(\mathbf{x_i}, \mathbf{\xi_i}(d_k) \right).
\end{equation}
Estimation proceeds as above using $\widehat{y^{T}_{G}}(d_k)$ in place of $\widehat{y^{T}}(d_k)$ to estimate $y^T(d_k)$. \citet{middleton_aronow11} and \citet{aronow_middleton_unbiased} 
demonstrate that $\widehat{y^{T}_{G}}(d_k)$ is location invariant.  Variance estimation proceeds as in Section \ref{varests}, using $\epsilon_i(d_k)$ in place of $y_i(d_k)$ so long as $g\left(\mathbf{x_i}, \mathbf{\xi_i}(d_k) \right)$ is fixed. 

An approximation to the difference estimator is given by regression adjustment using the data at hand. Regression can be thought of as a way to automate selection of the parameters in the difference estimator.  In doing so, unbiasedness is compromised although the regression estimator is typically consistent \citep[pp. 225-239]{sarndal_etal92}.  We may use weighted least squares to estimate a sensible parameter vector. For some common experimental designs, the least squares criterion will be optimal \citep{lin11_freedmans_critique}, and weighting by $1/\pi_i(d_k)$ ensures that the regression proceeds on a sample representative of the population of potential outcomes. With additional details on $\I_k$ and $g(.)$, it is possible to estimate optimal parameter vectors \citep[pp. 219-244]{sarndal_etal92}, though such values will typically be close to those produced by the weighted  least squares estimator (barring unusual and extreme forms of clustering).

Define an estimated parameter vector associated with exposure condition $d_k$
$$\mathbf{\widehat \xi}(d_k) = \arg\min\limits_{ \mathbf{\xi}(d_k)} \sum\limits_{i: D_i = d_k} \frac{1}{\pi_i(d_k)} \left[ Y_i - g\left(\mathbf{x_i}, \mathbf{\xi}(d_k) \right) \right]^2,$$
where $g(.)$ is the specification for the regression of $Y_i$ on  $\I(D_i=d_k)$ and $\mathbf{x_i}$. 
Then the regression estimator for the total is 
\begin{equation}
\label{eq:regression}
\widehat{y^{T}_{R}}(d_k) = \sum_{i=1}^N \I(D_i=d_k)\frac{Y_i - g\left(\mathbf{x_i}, \mathbf{\widehat \xi}(d_k) \right)}{\pi_i(d_k)}  +  \sum\limits_{i=1}^Ng\left(\mathbf{x_i}, \mathbf{\widehat \xi}(d_k) \right).
\end{equation}
Estimation proceeds as above using $\widehat{y^{T}_{R}}(d_k)$ in place of $\widehat{y^{T}_{HT}}(d_k)$ to estimate $y^T(d_k)$. Under weak regularity conditions on $g(.)$, a variance estimator based on a Taylor linearization of $\widehat{y^{T}_{R}}(d_k)$ is consistent \citep[pp.  236-237]{sarndal_etal92}.  The linearized variance estimator can be computed by substituting the residuals, $e_i = y_i(d_k) - g(\mathbf{x_i}, \mathbf{\widehat \xi}(d_k))$, for the $y_i(d_k)$ terms in constructing the variance estimator given in expression \eqref{eq:ate_var_estimator}. 

\subsection{Hajek Ratio Estimation Via Weighted Least Squares} 

The \citet{hajek71} ratio estimator is a refinement of the standard Horvitz-Thompson estimator that often facilitates efficiency gains at the cost of some finite $N$ bias and complications in variance estimation. Let us first consider the problem that the Hajek estimator is designed to resolve. The high variance of $\widehat{\mu_{HT}}(d_k)$ is often driven by the fact that some randomizations may yield units with exceptionally high values of the weights $1/\pi_i(d_k)$. The Hajek refinement allows the denominator of the estimator to vary according to the sum of the weights $1/\pi_i(d_k)$, thus shrinking the magnitude of the estimator when its value is large, and raising the magnitude of the estimator when its value is small. The Hajek ratio estimator is

\begin{equation}
\widehat{\mu_{H}}(d_k) = \frac{\sum_{i=1}^N \I(D_i=d_k)\frac{Y_i}{\pi_i(d_k)}}{\sum_{i=1}^N \I(D_i=d_k)\frac{1}{\pi_i(d_k)}}. \label{eq:hajek}
\end{equation}
Note that $\E[\sum_{i=1}^N \I(D_i=d_k)\frac{1}{\pi_i(d_k)}] = N$, so that the Hajek estimator is the ratio of two unbiased estimators. It is well known that the ratio of two unbiased estimators is not an unbiased estimator of the ratio. However, the bias will tend to be small relative to the estimator's sampling variability, and we may place bounds on its magnitude. 

By \citet{hartley54} and \citet[p. 176]{sarndal_etal92}, 
$$
\left|\E[\widehat{\mu_{H}}(d_k)]-\mu(d_k)\right|
\leq
\sqrt{\Var\left(\frac{1}{N}\sum_{i=1}^N \I(D_i=d_k)\frac{1}{\pi_i(d_k)}\right){\Var\left(\widehat{\mu_{H}}(d_k)\right)}}
$$
Under Conditions \ref{cond1} and \ref{cond2}, both variances will converge to zero, and the bias ratio will converge to zero. Practically speaking, the Hajek estimator can be computed using weighted least squares, with covariance adjustment through weighted least squares residualization.  Variance estimation proceeds via Taylor linearization \citep[pp. 172-176]{sarndal_etal92}.  A linearized variance estimator can be computed by substituting the residuals, $u_i = y_i(d_k) - \widehat{\mu_{H}}(d_k)$, for the $y_i(d_k)$ terms in constructing the variance estimator given in expression \eqref{eq:ate_var_estimator}. 

\section{Misspecification}

Recall Condition \ref{cond:sutva}, which states that the exposure mapping fully characterizes interference. Here we examine what happens when this assumption fails; e.g., there is interference between units that is not fully characterized by the exposure mapping. By ``misspecification'' of the exposure mapping, we refer to the situation in which the condition $D_i = d_k$ may be consistent with multiple potential outcomes for some $i$. As in Section \ref{avgpo}, we have randomization potential outcomes for unit $i$ as $y_i^r(\mathbf{z})$ for all $\mathbf{z} \in \Omega$.  

\begin{cond}[Misspecification]  There exists some $i \in \{1,...,N\}$ and $\z, \z' \in \Omega$ such that $f(\z,\theta_i) = f(\z',\theta_i)$ and $y_i^r(\z) \neq y_i^r(\z')$. Then $Y_i = \sum_{\z \in \Omega}  \mathbf{I}(\Z = \z) y_i^r(\z), \forall i \in \{ 1,...,N \}.$
\label{cond:misspec}
\end{cond}

 The following proposition shows the implications of misspecification for the potential outcome population total estimator given in expression \eqref{eq:ht_esimator}.

\begin{prop}\label{prop:misspec}
Define $\widehat{y^T_{HT}}(d_k)$ as above but suppose Condition \ref{cond:misspec} instead of Conditions \ref{cond:sutva} and \ref{cond:consistency}.  Then,
\begin{equation}
\E[\widehat{y^T_{HT}}(d_k)]  = \sum_{i=1}^N \sum_{\mathbf{z}:f(\mathbf{z},\theta_i)=d_k} w_{i,\z} y_i^r(\mathbf{z}) \label{eq:misspec},
\end{equation}  
where $ w_{i,\z} = {p_{\mathbf{z}}}/ {\pi_i(d_k)}$.
\end{prop}
Under Condition \ref{cond:misspec}, the estimator $\widehat{\mu_{HT}}(d_k) = (1/N)\widehat{y^T_{HT}}(d_k)$ is unbiased for the population mean of what Hudgens and Halloran (2008, p. 834) refer to as the ``individual average potential outcome'' given $D_i = d_k$.  The causal effect estimate given in \eqref{eq:ht_causal_effect}, which compares mean outcomes given exposures $d_k$ versus $d_l$, is a difference in population means of individual average randomization potential outcomes given different restrictions on the set of treatments implied in constructing exposures $d_k$ and $d_l$.

\begin{corollary}
Under Condition \ref{cond:misspec},
$$
\E[\widehat{\tau_{HT}}(d_k, d_l)] = \frac{1}{N} \sum_{i=1}^N  \left[\sum_{\mathbf{z}:f(\mathbf{z},\theta_i)=d_k}   w_{i,\z} y_i^r(\mathbf{z})  - \sum_{\mathbf{z'}:f(\mathbf{z'},\theta_i)=d_l} w_{i,\z'} y_i^r(\mathbf{z}') \right].
$$
\end{corollary}
Inference for such an effect would not follow immediately from the results above. However, under partial interference, inference would follow from the results of \citet{liu-hudgens2014-interference-intervals}.

\section{A naturalistic simulation with social network data}\label{sec:simulation}

We use a naturalistic simulation to illustrate how our framework may be applied and also to study operating characteristics of the proposed estimators in a realistic sample. We estimate direct and indirect effects of an experiment with individuals linked in a complex, undirected social network.  We use friendship network data from American school classes collected through the National Longitudinal Study of Adolescent Health (Add Health). 
 The richness of these data makes Add Health a canonical dataset for methodological research related to social networks, as with  
\citet{bramoulle_etal2009_peer_effects},
\citet{chung_etal2008_latent_transition_analysis},
\citet{goel_salganik10_assessing_rds}, 
\citet{goodreau_etal2009_ergm},
\citet{goodreau2007_ergms},
\citet{handcock_etal2007_network_clustering}, and
\citet{hunter_etal2008_social_network_models_fit}.
We simulate experiments in which a treatment, $\Z$, is randomly assigned without replacement and with uniform probability to $1/10$ of individuals in a school network.  Indirect effects are transmitted only within a subject's school.  This simulated experiment resembles various studies of network persuasion campaigns \citep{chen_etal2010_diffusion, aral_walker2011_viral,paluck2011_peer_pressure}, 
including the field experiment that we analyze below. 

We define the exposure mapping as a function $f(\z, \theta_i)$ such that the parameter, $\theta_i$, is a column vector equal to the transpose of subject $i$'s row in a network adjacency matrix (modified such that we have zeroes on the diagonal).  The inner product, $\z'\theta_i$, counts the number of subject $i$'s peers assigned to treatment.  We use a simple exposure mapping that captures direct and indirect effects of the treatment, with indirect effects being transmitted to a subject's immediate peers:
\begin{align}
f(\z, \theta_i) = \left\{\begin{array}{lr} 
  d_{11}  \text{ }( \text{Direct + Indirect Exposure}):& z_i\I(\z'\theta_i>0)  = 1 \\  
  d_{10}  \text{ }( \text{Isolated Direct Exposure}): &  z_i \I(\z'\theta_i=0) = 1 \\ 
   d_{01}   \text{ }( \text{Indirect Exposure}):&  (1-z_i)\I(\z'\theta_i>0) = 1  \\  
   d_{00}  \text{ }( \text{No Exposure}): &(1-z_i)\I(\z'\theta_i=0)  = 1 \end{array} \right. \nonumber
\end{align}
where each unit falls into exactly one of the four exposure conditions. This exposure mapping was selected to mimic the one used in the application studied in the next section.  A contrast of mean outcomes under $d_{10}$ versus $d_{00}$ isolates the effect of direct exposure in the absence of any interaction with indirect exposure, whereas the $d_{11}$-$d_{00}$ contrast yields an effect that incorporates such interactive effects.  The $d_{01}$-$d_{00}$ contrast isolates the effect of indirect exposure in the absence of any interaction with direct exposure.

This experiment is repeated independently across the 144 school classes included in Add Health, with an average class size of 626 students. We constructed the school network graphs as undirected graphs where a link between two students was assigned if either student nominated the other as a friend in the AddHealth survey.  Students could nominate up to 5 male and 5 female friends.  To ensure that our effect estimates all refer to the same underlying population, we dropped subjects that reported zero friendship ties.  For the resulting sample, 42\% of students have network degree in the 1 to 5 range, 40\% in the 6 to 10 range, 18\% in the 11 to 20 range, and 1\% greater than 20, with a maximum degree of 39. To give an idea of the range of exposure probabilities, for the student with degree of 39, the probability of isolated direct exposure was 0.00067.  In the appendix, we display the cumulative distribution functions for the four exposure probabilities.  About 3\% of students have an exposure probability of less than 0.01 for the direct + indirect exposure condition, 0.5\% for isolated direct exposure, and then there were no cases with probabilities less than 0.01 for either the indirect- or no-exposure conditions.

\begin{figure}[!t]
\begin{center}
\includegraphics[width=1\textwidth]{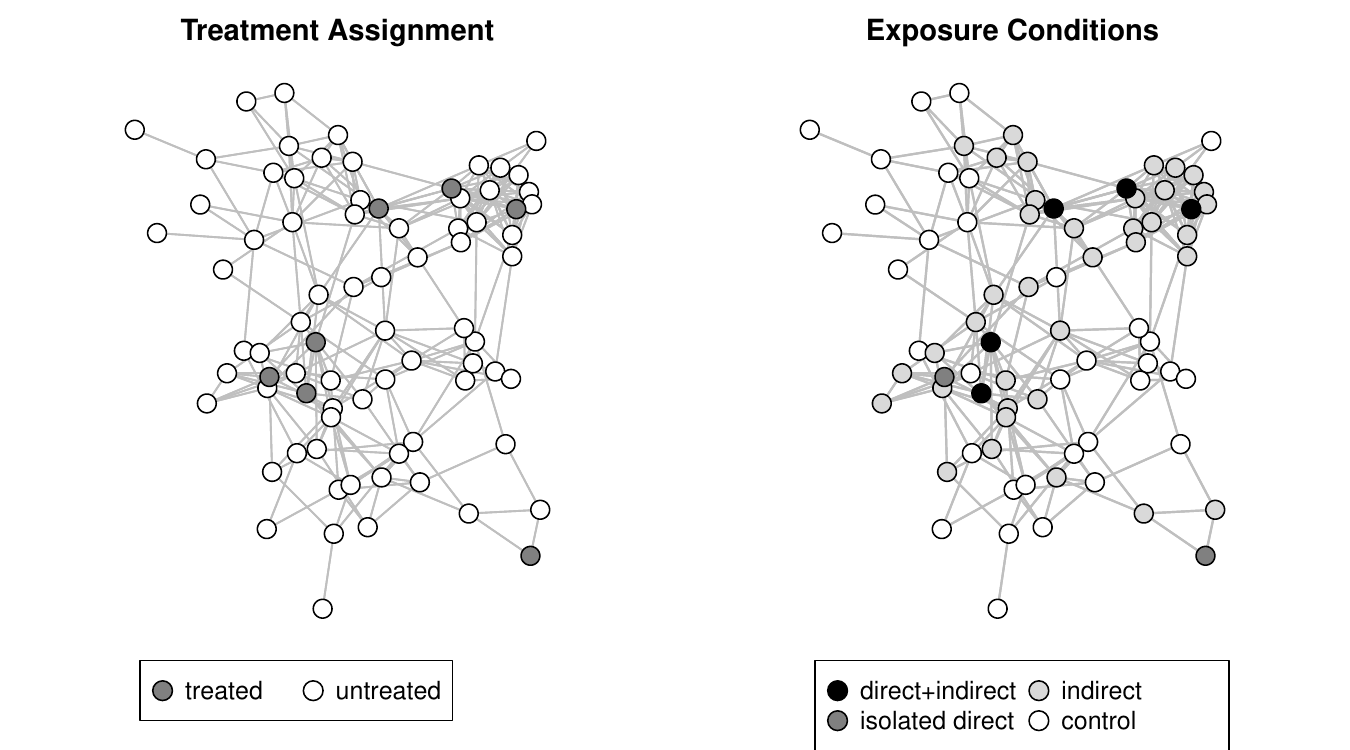}
\end{center}
\vspace{-2em}
\caption{Illustration of a treatment assignment (left) and then treatment-induced exposures (right) for one of the school classes in the study. Each dot is a student, and each line represents an undirected friendship tie.\label{fig:network}}
\end{figure}

Figure \ref{fig:network} illustrates a treatment assignment and corresponding treatment-induced exposures under this mapping.  The figure illustrates two key issues that our methods address.  First is the connection between a unit's underlying traits, in this case its network degree, and propensity to fall into one or another exposure condition.  The second is the irregular clustering that occurs in exposure conditions.  Such irregular clustering is precisely what one must address in deriving variance estimates and intervals for estimated effects.  

We use as our outcome a variable in the dataset that records the number of after-school activities in which each student participates.  This variable defines the $y_i(d_{00})$ values---that is, potential outcomes under the ``control'' exposure.  This makes our simulation naturalistic not only in the networks that define the interference patterns, but also in the outcome data.  The variable exhibits a high degree of right skew, with mean 2.14, standard deviation 2.64, and 0, .25, .5, .75, and 1 quantiles of 0, 1, 2, 3, and 33, respectively. 
We consider a simple ``dilated effects'' scenario \citep{rosenbaum1999_dilated_effects_sensitivity} where potential outcomes are such that $y_i(d_{11})= 2\times y_i(d_{00}), y_i(d_{10})= 1.5 \times y_i(d_{00}), y_i(d_{01})= 1.25\times y_i(d_{00})$.  We run 500 simulated replications of the experiment, applying five estimators in each scenario:
\begin{itemize}
\item The Horvitz-Thompson estimator for the causal effect given in expression \eqref{eq:ht_causal_effect}, with the associated conservative variance estimator, given in expression \eqref{eq:ate_var_estimator};
\item The Hajek ratio estimator given in expression \eqref{eq:hajek}, with the associated linearized variance estimator;
\item The weighted least squares (WLS) estimator given in expression \eqref{eq:regression}, adjusting for network degree as the sole covariate, with the associated linearized variance estimator;
\item An ordinary least squares (OLS) estimator that regresses the outcome on indicator variables for the exposure conditions, adjusting for network degree as a covariate, with \citet{mackinnon_white1985_hc2}'s finite sample adjusted ``HC2'' heteroskedasticity consistent variance estimator;
\item A simple difference in sample means (DSM) for the exposure conditions, also with the HC2 estimator.
\end{itemize}
With respect to point estimates, the Horvitz-Thompson estimator is unbiased but possibly unstable, while the Hajek and WLS estimators are consistent and expected to be more stable.  The DSM estimator is expected to be biased because it totally ignores relationships between exposure probabilities and outcomes.  The OLS estimator controls for network degree, and so this will remove bias due to correlation between exposure probabilities and outcomes. However, OLS is known to be biased in its aggregation of unit-level heterogeneity in causal effects \citep{angrist_krueger99_handbook}.  With respect to standard error estimates and confidence intervals, the variance estimators for the Horvitz-Thompson, Hajek, and WLS estimators are expected to be conservative though informative.  The variance estimators for OLS and DSM may be anti-conservative because they ignore the clustering in exposure conditions.


Table \ref{tab:simresults_dilated} shows results of the simulation study, which conform to expectations.  The Horvitz-Thompson, Hajek, and WLS estimators exhibit no perceivable bias.  The Horvitz-Thompson estimator exhibits higher variability than the Hajek and WLS estimators, although the differences are not very pronounced, perhaps owing to the small number of cases with very small exposure probabilities.  The OLS estimator and DSM estimator are heavily biased when considered relative to the variability of the effect estimates.  The bias in OLS is expected because unit-level causal effects, defined in terms of differences, are heterogenous from unit to unit when underlying potential outcomes are based on dilated effects.  Thus OLS will suffer from an aggregation bias in addition to any biases due to inadequate conditioning on network degree.  The standard error estimates for the Horvitz-Thompson, Hajek, and WLS estimators are informative but conservative, resulting in empirical coverage rates that exceed nominal levels.  The intervals for the OLS and DSM variance estimators badly undercover, primarily due to the bias in the point estimates rather than understatement of variability.

\begin{table}[!t]
\caption{Results from school friends' network simulated experiment}
\label{tab:simresults_dilated}
\begin{center}

\begin{tabular}{llrrrrrrr}
\hline
          &                   &     &       &       &  Mean &		  95\% CI  &  90\% CI  \\
Estimator  &Estimand             &Bias&  S.D. &  RMSE &  S.E. & Coverage &     Coverage\\
\hline
HT   & $\tau(d_{01},d_{00})$  & {0.00 }&  0.04 &  0.04 &  0.05 &  0.960    & 0.924\\
 & $\tau(d_{10},d_{00})$      & {0.00 }&  0.10 &  0.10 &  0.19 & 0.986    & 0.970\\
 & $\tau(d_{11},d_{00})$      & {0.00 }&  0.13 &  0.13 &  0.28 &  0.990    &  0.970\\
\hline
Hajek & $\tau(d_{01},d_{00})$ & {0.00 }&  0.03 &  0.03 &  0.03 &  0.968    & 0.916\\
 & $\tau(d_{10},d_{00})$      & {0.00 }&  0.07 &  0.07 &  0.13 &  0.992    & 0.970\\
 & $\tau(d_{11},d_{00})$      & {0.00 }&  0.12 &  0.12 &  0.25 & 0.986    & 0.970\\ 
\hline
WLS & $\tau(d_{01},d_{00})$   & {0.00 }&  0.03 &   {0.03} &  0.03 &  0.970    &  0.928\\
 & $\tau(d_{10},d_{00})$      & {0.00 }&  0.07 &   {0.07} &  0.12 &  0.992    & 0.968\\
 & $\tau(d_{11},d_{00})$      & {0.00 }&  0.11 &   {0.11} &  0.25 &  0.988    & 0.950\\
\hline
OLS & $\tau(d_{01},d_{00})$   &  {-0.02}&  0.03 &  0.03 &  0.02 & {0.842}    & {0.768}\\
 & $\tau(d_{10},d_{00})$      &  {-0.08}&  0.06 &  0.10 &  0.07 &   {0.706}    &   {0.576}\\
 & $\tau(d_{11},d_{00})$      &  {0.12} &  0.09 &  0.15 &  0.09 &   {0.660}    & {0.530}\\
\hline
DSM & $\tau(d_{01},d_{00})$    &  {0.42} &  0.02 &  0.42 &  0.02 &   {0.000}      &  {0.000}\\
 & $\tau(d_{10},d_{00})$      &  {-0.08}&  0.06 &  0.10 &  0.07 &   {0.726}    &   {0.614}\\
 & $\tau(d_{11},d_{00})$      &  {0.56} &  0.09 &  0.57 &  0.09 &   {0.000}    &   {0.000}\\
 \hline
\end{tabular}
\end{center}
\begin{scriptsize}
HT = Horvitz-Thompson estimator with conservative variance estimator.\\
Hajek = Hajek estimator with linearized variance estimator.\\
WLS = Least squares weighted by exposure probabilities with covariate adjustment for network degree and linearized variance estimator.\\
OLS = Ordinary least squares with covariate adjustment for network degree and heteroskedasticity consistent variance estimator.\\
DSM = Simple difference in sample means with no covariate adjustment and heteroskedasticity consistent variance estimator.\\
S.D. = Empirical standard deviation from simulation; RMSE = Root mean square error; S.E. = standard error estimate; CI = Normal approximation confidence interval.
\end{scriptsize}
\end{table}

\section{Analysis of a social network field experiment}\label{sec:experiment}

\begin{figure}[!ht]
\begin{center}
\includegraphics[width=.8\textwidth]{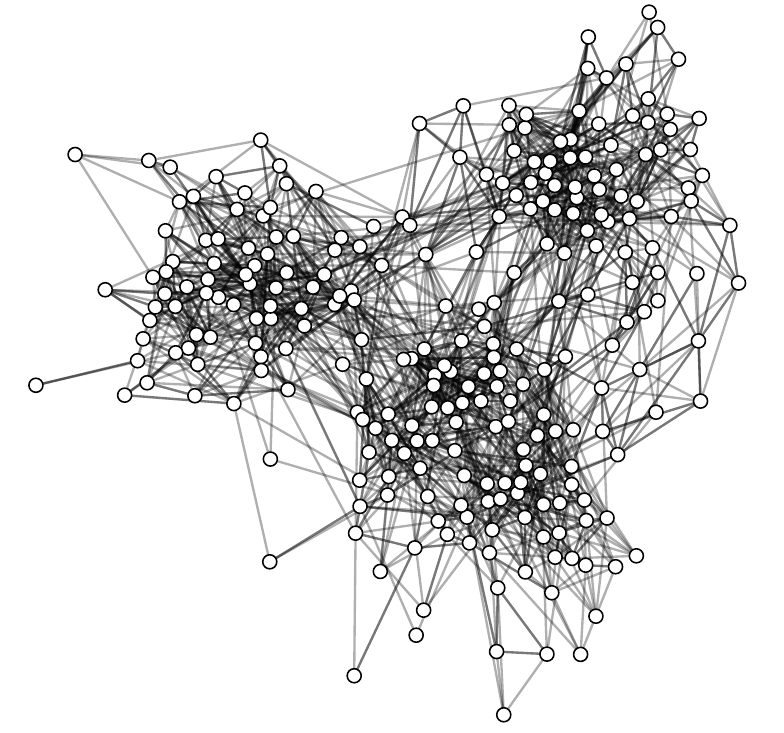}
\caption{Example of a proximity network for one school in the social network field experiment. Network edges were measured using student nomination data in the first survey.}
\label{samplenetwork}
\end{center}
\end{figure}

In this section we analyze a field experiment on the promotion of anti-conflict norms and behavior among American middle school students.  The experiment sought to shed light on how such a program might affect attitudes and behaviors of participant youth and also, crucially, to understand how these effects transmit to participants' social network peers. Full details and a richer analysis of the experiment are given in \citet{paluck-etal2016-pnas-highschool}.  The experiment involved two levels of randomization. First, 28 of 56 schools were randomly selected to host the anti-conflict program, via block randomization.  Within all schools, a group of between 40 to 64 students were nonrandomly selected as eligible to participate in the program.  Within each school hosting the program, half of the eligible students were then block randomized to participate in the program, with blocking on gender, grade, and a measure of network closure.  Every two weeks over the course of the school year, the program had participants attend meetings with program staff during which they discussed social conflicts and patterns of exclusion at their school and formulated behavioral strategies to help friends and other students.  At the beginning of the school year, the research team measured students' social networks, asking students to nominate up to 10 students in their school that they had chosen to spend time with, face to face or online, in the last two weeks. These nominations were used to construct an undirected adjacency graph, so that students were considered ``peers'' if either student nominated one another. In Figure \ref{samplenetwork}, we present an illustrative graph of one of these networks. As expected, students of the same grade and gender are more likely to associate with one another. At the end of the school year, the research team implemented a survey to measure behaviors and attitudes that reflected conflict-related norms.  In the current analysis, we focus on one particular behavior: (self-reports of) wearing a wristband issued to students through the program that was meant to reflect a student's public endorsement of anti-conflict norms.  

Given this design, let $z_i=0,1$ be an indicator for whether student $i$ is assigned to participate in the program and let $\z$ be the vector of student-level assignments in the school.  Let $s_i=0,1$ be an indicator for whether subject $i$'s school hosts the program. Finally, as in the simulation study above, let $\theta_i$ be a column vector equal to the transpose of student $i$'s row in the school network adjacency matrix (with zeroes on the diagonal), in which case $\z'\theta_i$, is again the number of subject $i$'s peers that are assigned to participate in the program.  Then, we define the exposure mapping as follows:
\begin{align}
f(\z, \theta_i) = \left\{\begin{array}{lr} 
d_{111}  \text{ }( \text{Direct + Indirect Exposure}): & z_i\I(\z'\theta_i>0)s_i = 1 \\
d_{101}   \text{ }( \text{Isolated Direct Exposure}):& z_i\I(\z'\theta_i=0)s_i  = 1  \\  
d_{011}  \text{ }( \text{Indirect Exposure}): &(1-z_i)\I(\z'\theta_i >0) s_i  = 1 \\  
d_{001}  \text{ }( \text{School Exposure}):& (1-z_i)\I(\z'\theta_i=0)s_i  = 1 \\  
d_{000} \text{ }( \text{No Exposure}):  & (1-s_i) = 1
\end{array} \right. \nonumber
\end{align}
Three features emerge from examination of the exposure mapping. First, the exposure mapping reflects four different sources of exposure to the program: being in a school with the program (School), having a peer who was a participation student (Indirect), and being a participating student (Direct). Second, only students selected as ``eligible'' have a non-zero possibility of being in all exposure conditions.  Thus, we limit the present analysis to the set of students with nonzero probabilities of exposure ($N=2,050$). (\citet{paluck-etal2016-pnas-highschool} examine effects for members of the ineligible subpopulation who nonetheless have non-zero probability of indirect exposure.) Third, the conditions for our asymptotic results hold in the number of schools.  This exposure model provides a parsimonious characterization of first-order peer effects and school-wide climate effects, which were the primary effects of interest for \citet{paluck-etal2016-pnas-highschool} when designing the experiment.  If other types of peer effects are present, analysis under this exposure model estimates treatment-regime specific aggregates that average over those other effects, as described in the section on misspecification above.

\begin{table}[!t]
\caption{Social network experiment results: effects of exposures on probability of wearing a program wristband}
\label{tab:expresults}
\begin{center}
\begin{tabular}{llrrc}
\hline
Estimator  &Estimand          &Estimate&  S.E. &  95\% CI\\
\hline
HT   & $\tau(d_{001},d_{000})$  & 0.057 & 0.062 &  (-0.065, 0.179) \\
 & $\tau(d_{011},d_{000})$      & 0.154& 0.029 &  (0.097, 0.211) \\
 & $\tau(d_{101},d_{000})$      & 0.305 & 0.141 &  (0.029, 0.581) \\
 & $\tau(d_{111},d_{000})$      & 0.299 & 0.020 &  (0.260, 0.338) \\
\hline
Hajek   & $\tau(d_{001},d_{000})$  & 0.058 & 0.064 &  (-0.067, 0.183) \\
 & $\tau(d_{011},d_{000})$      & 0.154 & 0.037 &  (0.081, 0.227) \\
 & $\tau(d_{101},d_{000})$      & 0.292 & 0.123 &  (0.051, 0.533) \\
 & $\tau(d_{111},d_{000})$      & 0.307 & 0.049 &  (0.211, 0.403) \\
\hline
WLS   & $\tau(d_{001},d_{000})$  & 0.056 & 0.066 &  (-0.072, 0.186) \\
 & $\tau(d_{011},d_{000})$      & 0.156 & 0.037 &  (0.083, 0.229) \\
 & $\tau(d_{101},d_{000})$      & 0.295 & 0.124 &  (0.050, 0.536) \\
 & $\tau(d_{111},d_{000})$      & 0.306 & 0.049 &  (0.212, 0.404) \\
\hline
\end{tabular}
\end{center}
\begin{scriptsize}
HT = Horvitz-Thompson estimator with conservative variance estimator.\\
Hajek = Hajek estimator with linearized variance estimator.\\
WLS = Least squares weighted by exposure probabilities with covariate adjustment for network degree and linearized variance estimator.\\
S.E. = Estimated standard error; CI = Normal approximation confidence interval.
\end{scriptsize}
\end{table}

Table \ref{tab:expresults} presents Horvitz-Thompson (HT), Hajek, and weighted least squares (WLS) estimates of the effects of different exposure conditions relative to the no exposure condition.  The WLS estimates control for a subject's network degree (as in the simulation study above).  The outcome of interest, $y_i\in \{0,1\}$ is a binary indicator for whether the subject wore a program wristband. The effect estimates characterize, for eligible students, the average increase in the probability of wearing a wristband relative to the average in the no exposure condition. (The average for eligible students in the no exposure condition was essentially zero, at 0.000.)

The HT, Hajek, and WLS results all mostly agree.  They suggest that being in a program school but being a non-participant with no participant peers ($d_{001}$) has negligible effects for eligible students: our point estimate suggests about a six percentage point increase in the probability of wearing a wrist band, although the 95\% confidence interval has a lower bound of about -7 percentage points.  However, effects for eligible students with either indirect or direct exposure are substantially larger.  The effect of indirect exposure ($d_{011}$) is about a 15 to 16 percentage point increase in the probability of wearing a wrist band (95\% confidence interval between about 8 and 23 percentage points). The effect of direct exposure, whether or not it is accompanied by indirect exposure ($d_{101}$ or $d_{111}$) is about a 30 percentage point increase in the probability of wearing a wrist band (95\% confidence interval between about 5 and 50 percentage points for the $d_{101}$ condition and about 21 and 40 percentage points for the $d_{111}$ condition).  

Thus, the program is seen as having substantial direct but also indirect effects on subject's willingness to endorse anti-conflict norms by wearing a program wristband.  These indirect effects mean that one would drastically underestimate the effect of the program if one performed a naive analysis that simply compared participant and non-participant individuals in schools hosting the program.  Moreover, an analysis that failed to account for indirect effects might understate the cost-effectiveness of the program: substantial increases in school-level expressions of commitment to anti-conflict norms would not require administering the program to everyone.

\section{Conclusion} \label{interfereconc}

This paper proposes an analytical framework for causal inference under interference
and applies it to the analysis of experiments on social networks.  As discussed in the introduction, the framework can be applied to other settings where interference is considered to be important.
The framework integrates (i) an experimental design that defines the probability distribution for treatments assigned, (ii) an exposure mapping that relates treatments assigned to exposures received, and (iii) an estimand chosen to make use of an experimental design to answer questions of substantive interest.  Using this framework, we develop methods 
for estimating average unit-level causal effects 
of exposures from a randomized experiment.
Our approach combines the known randomization process with the analyst's definition of treatment exposure, thus permitting inference under clear and defensible assumptions. Importantly, the union of the design of the experiment and the exposure mapping may imply unequal probabilities of exposure and forms of dependence between units that may not be obvious ex ante.

We develop estimators based on results from the literature on unequal probability sampling rooted in the foundational insights of \citet{horvitz_thompson}.  The estimators are derived from the known sampling distribution of the ``direct'' treatment, and they provide a basis for unbiased effect estimation and conservative variance estimation.  Wald-type intervals based on a normal approximation provide a reasonable reflection of large $N$ behavior when clustering of exposure indicator values is limited.  Nonetheless, it is well known that Horvitz-Thompson-type estimators may be volatile in cases where selection probabilities vary greatly or exhibit strong inverse correlation with outcome values \citep{basu1971_elephants}.   Thus, we provide refinements that allow for variance control via covariance adjustment and Hajek estimation. 

Our approach combines 
minimal assumptions about restrictions on potential outcomes
with randomization-based estimators, and may be characterized as design-consistent.  
The estimands and methods presented here may be useful in evaluating alternative experimental designs for estimating causal effects in the presence of interference \citep{airoldi2016-interference, baird-etal-spillover-design, eckles-etal2014, toulis-kao2013, ugander-etal2013}.
Finally, the framework is readily applicable to deriving estimators for estimands other than the average unit-level effect of exposures.

\appendix

\section{Proofs}

\subsection{Proof of Proposition \ref{prop:sim-converge}}
The replication procedure is equivalent to drawing a random sample without replacement from $\Omega$, with probabilities of selection equal to those which are defined in the randomization plan. The result follows from the Strong Law of Large Numbers.

\subsection{Proof of Lemma \ref{lemma:ht}}
To show unbiasedness, by \ref{cond:consistency} we have
\begin{align}
\E[\widehat{y^T_{HT}}(d_k)] & = \E\Big[\sum_{i=1}^N \I(D_i=d_k)\frac{Y_i}{\pi_i(d_k)} \Big] \nonumber \\ & = \sum_{i=1}^N \E[\I(D_i=d_k)]\frac{y_i(d_k)}{\pi_i(d_k)} = \sum_{i=1}^N y_i(d_k). \nonumber
\end{align}
The variance expression follows from the fact that $\widehat{y^T_{HT}}(d_k)$ is a sum of correlated random variables. 

\subsection{Proof of Proposition \ref{prop:bias-bound}}
We have,
\begin{align}
\E[\widehat{y^T_{HT,R}} (d_k)] & = \sum_{i=1}^N y_i(d_k)\pi_i(d_k)  \E\left[ \frac{R+1}{X_i+1}\right] \nonumber \\ 
& = \sum_{i=1}^N y_i(d_k) \left[1-(1-\pi_i(d_k))^{R+1} \right] \nonumber \\ 
& = y^T(d_k) - \sum_{i=1}^N y_i(d_k)(1-\pi_i(d_k))^{R+1}. \nonumber
\end{align}
And so,
\begin{align}
\left|\E[\widehat{y^T_{HT,R}} (d_k)] - y^T(d_k) \right| & = \left| \sum_{i=1}^N y_i(d_k)(1-\pi_i(d_k))^{R+1}\right| \nonumber 
\\
& \le | y^T(d_k) | (1-\pi_0(d_k))^{R+1}.\nonumber
\end{align}

\subsection{Proof of Proposition \ref{prop:ht}}
Results \eqref{eq:etauhat} and \eqref{eq:tru_var} follow from Lemma \ref{lemma:ht} and properties of the variance operator. For the covariance term \eqref{eq:totals_covariance}, first note that $\pi_{ii}(d_k,d_l) = 0$. Then following \citet{wood08},
\begin{align}
\Cov[\widehat{y^T_{HT}}(d_k),\widehat{y^T_{HT}}(d_l)] & =  \sum_{i=1}^N \sum_{j=1}^N \Cov \left[\I(D_i=d_k),\I(D_j=d_l)\right]\frac{y_i(d_k)}{\pi_i(d_k)} \frac{y_j(d_l)}{\pi_j(d_l)} \nonumber\\
& = \sum_{i=1}^N \sum_{j=1}^N  \frac{y_i(d_k)}{\pi_i(d_k)} \frac{y_j(d_l)}{\pi_j(d_l)} \left[\pi_{ij}(d_k,d_l)- \pi_i(d_k)\pi_j(d_l) \right].\nonumber
\end{align}

\subsection{Key results for Propositions \ref{lem:varbias}, \ref{a2}, \ref{ncov}, and \ref{prop:no-bias-cov}}

We reproduce key results from Aronow and Samii (2013) for the conservative variance corrections.  We do so in the general case of the Horvitz-Thompson estimator for a population total.  Suppose a population $U$ indexed by $1,...,k,...,N$ and a sampling design such that the probability of inclusion in the sample for unit $k$ is given by $\pi_k$, and the joint inclusion probability for units $k$ and $l$ is given by $\pi_{kl}$.

The Horvitz-Thompson estimator of a population total is given by $$\hat{t} = \sum_{k \in s} \frac{y_k}{\pi_k} = \sum_{k \in U} I_k \frac{y_k}{\pi_k},$$ where $I_k \in \{0,1\}$ is unit $k$'s inclusion indicator, the only stochastic component of the expression, with $\text{E}(I_k) = \pi_k$, the inclusion probability, and $s$ and $U$ refer to the sample and the population, respectively.  Define $\text{E}(I_kI_l)= \pi_{kl}$, the probability that both units $k$ and $l$ from $U$ are included in the sample.  Since $I_kI_k= I_k$, $\text{E}(I_kI_k)=\pi_{kk}=\pi_k$ by construction. The variance of the Horvitz-Thompson estimator for the total is given by
\begin{align}
\text{Var}(\hat{t}) & = \sum_{k \in U} \sum_{l \in U} \text{Cov}(I_k,I_l) \frac{y_k}{\pi_k} \frac{y_l}{\pi_l} \nonumber \\
& = \sum_{k \in U}\text{Var}(I_k)\left(\frac{y_k}{\pi_k}\right)^2 + \sum_{k \in U} \sum_{l \in U \backslash k }  \text{Cov}(I_k,I_l) \frac{y_k}{\pi_k} \frac{y_l}{\pi_l}. \nonumber 
\end{align}
Under a measurable design, two conditions obtain: (1) $\pi_k>0$ and $\pi_k$ is known for all $k \in U$ and (2) $\pi_{kl} >0$ and $\pi_{kl}$ is known for all $k,l \in U$.  Non-measurable designs include those for which either of the two conditions for a measurable design do not hold. We label a sample from a measurable design, $s^M$, and an unbiased estimator for $\text{Var}(\hat{t})$ on $s^M$ is given by,
$$
\widehat{\text{Var}}(\hat{t}) = \sum_{k \in s^M} \sum_{l \in s^M} \frac{\text{Cov}(I_k,I_l)}{\pi_{kl}}\frac{y_k}{\pi_k} \frac{y_l}{\pi_l} = \sum_{k \in U} \sum_{l \in U} I_k I_l \frac{\text{Cov}(I_k,I_l)}{\pi_{kl}}\frac{y_k}{\pi_k} \frac{y_l}{\pi_l},
$$
where the only stochastic part of the latter expression is $I_k I_l$, and unbiasedness is due to $\text{E}(I_k I_l) =\pi_{kl}$.

Suppose a non-measurable design for which $\pi_{kl} = 0$ for some units $k, l \in U$. We label a sample from such a non-measurable design as $s^0$. Because $I_k$ is a Bernoulli random variable with probability $\pi_k$, $\text{Cov}(I_k,I_l) = \pi_{kl} - \pi_k\pi_l$ for $k\ne l$, and $\text{Cov}(I_k,I_k) = \text{Var}(I_k) = \pi_k(1-\pi_k)$. Then, we can re-express the variance above as,
\begin{align}
\text{Var}(\hat{t}) & = \sum_{k \in U}\pi_k(1-\pi_k) \left(\frac{y_k}{\pi_k}\right)^2 + \sum_{k \in U} \sum_{l \in U \backslash k }  (\pi_{kl} - \pi_k\pi_l) \frac{y_k}{\pi_k} \frac{y_l}{\pi_l} \nonumber \\
& = \sum_{k \in U}\pi_k(1-\pi_k) \left(\frac{y_k}{\pi_k}\right)^2 + \sum_{k \in U} \sum_{l \in \{U \backslash k : \pi_{kl}>0\} }  (\pi_{kl} - \pi_k\pi_l) \frac{y_k}{\pi_k} \frac{y_l}{\pi_l} - \underbrace{\sum_{k \in U} \sum_{l \in \{U \backslash k : \pi_{kl}=0\} } y_k y_l}_A . \nonumber
\end{align}
For $k$ and $l$ such that $\pi_{kl}=0$, the sampling design will never provide information on the component of the variance labeled as $A$ above, since we will never observe $y_k$ and $y_l$ together.   

When $\widehat{\text{Var}}(\hat{t})$ is applied to $s^0$, the result is unbiased for $\text{Var}(\hat{t}) + A$.  We state this formally as follows:
\begin{prop}[Aronow and Samii, 2013, Prop. 1]\label{prop:as-prop1}
When $s^0$ refers to a sample from a design with some $\pi_{kl}=0$, we have,
$$
\text{E}\left[\widehat{\text{Var}}(\hat{t})\right]= \text{Var}(\hat{t}) + \sum_{k \in U} \sum_{l \in \{U \backslash k : \pi_{kl}=0\} } y_k y_l = \text{Var}(\hat{t}) + A.
$$
\end{prop}

\begin{proof}
The result follows from,
\begin{align}
\text{E}\left[ \sum_{k \in s^0} \sum_{l \in s^0} \frac{\text{Cov}(I_k,I_l)}{\pi_{kl}}\frac{y_k}{\pi_k} \frac{y_l}{\pi_l} \right] & = \text{E} \left[ \sum_{k \in U} \sum_{l \in \{U: \pi_{kl}>0\}} I_k I_l \frac{\text{Cov}(I_k,I_l)}{\pi_{kl}}\frac{y_k}{\pi_k} \frac{y_l}{\pi_l} \right]  \nonumber \\
& = \sum_{k \in U} \text{Var}(I_k)\left(\frac{y_k}{\pi_k}\right)^2  + \sum_{k \in U} \sum_{l \in \{U\backslash k: \pi_{kl}>0\}}  \text{Cov}(I_k,I_l) \frac{y_k}{\pi_k} \frac{y_l}{\pi_l} \nonumber \\
& = \text{Var}(\hat{t}) + \sum_{k \in U} \sum_{l \in \{U \backslash k : \pi_{kl}=0\} } y_k y_l \nonumber \\
& = \text{Var}(\hat{t}) + A. \nonumber 
\end{align}
\end{proof}

Now, consider the following variance estimator:
$$
\widehat{\text{Var}_C}(\hat{t}) =  \sum_{k \in U} \sum_{l \in \{U: \pi_{kl}>0\}} I_k I_l \frac{\text{Cov}(I_k,I_l)}{\pi_{kl}}\frac{y_k}{\pi_k} \frac{y_l}{\pi_l} +  \sum_{k \in U} \sum_{l \in \{U \backslash k : \pi_{kl}=0\} } \left(I_k \frac{|y_k|^{a_{kl}}}{a_{kl} \pi_k} + I_l \frac{|y_l|^{b_{kl}}}{b_{kl} \pi_l}\right),
$$
where $a_{kl}, b_{kl}$ are positive real numbers such that $\frac{1}{a_{kl}} + \frac{1}{b_{kl}} = 1$ for all pairs ${k},{l}$ with $\pi_{kl} = 0$.  

\begin{prop}[Aronow and Samii, 2013, Prop. 2] \label{prop:as-prop2}
$$
\text{E} \left[\widehat{\text{Var}_C}(\hat{t}) \right] \geq \text{Var}(\hat{t}).
$$
\end{prop}

\noindent {\it Proof.} By Young's inequality,
$$
\frac{|y_k|^{a_{kl}}}{a_{kl}} + \frac{|y_l|^{b_{kl}}}{b_{kl}} \geq |y_k| |y_l|,
$$
if $\frac{1}{a_{kl}} + \frac{1}{b_{kl}} = 1$. Define $A^*$ such that,
$$
A^* =  \sum_{k \in U} \sum_{l \in \{U \backslash k : \pi_{kl}=0\} }  \frac{|y_k|^{a_{kl}}}{a_{kl}} + \frac{|y_l|^{b_{kl}}}{b_{kl}} \geq  \sum_{k \in U} \sum_{l \in \{U \backslash k : \pi_{kl}=0\} } |y_k| |y_l| \geq \sum_{k \in U} \sum_{l \in \{U \backslash k : \pi_{kl}=0\} } y_k y_l = A
$$
and
$$
A^* \geq  \sum_{k \in U} \sum_{l \in \{U \backslash k : \pi_{kl}=0\} } |y_k| |y_l| \geq \sum_{k \in U} \sum_{l \in \{U \backslash k : \pi_{kl}=0\} } -y_k y_l = -A.
$$
Therefore
$$
 \text{Var}(\hat{t}) + A  + A^* \geq \text{Var}(\hat{t}).
$$

The associated Horvitz-Thompson estimator of $A^*$ would be
$$
\widehat{A^*} = \sum_{k \in U} \sum_{l \in \{U \backslash k : \pi_{kl}=0\} } \left(I_k \frac{|y_k|^{a_{kl}}}{a_{kl} \pi_k} + I_l \frac{|y_l|^{b_{kl}}}{b_{kl} \pi_l}\right),
$$
which is unbiased by $\text{E}(I_k) = \pi_k$ and $\text{E}(I_l) = \pi_l$.

Since $\text{E}\left[\widehat{A^*}\right] = A^*$, by Proposition \ref{prop:as-prop1},
$$
\text{E}\left[ \sum_{k \in s^0} \sum_{l \in s^0} \frac{\text{Cov}(I_k,I_l)}{\pi_{kl}}\frac{y_k}{\pi_k} \frac{y_l}{\pi_l} + \widehat{A^*} \right] =  \text{Var}(\hat{t}) + A + A^*
$$
$$
\text{E} \left[\widehat{\text{Var}_C}(\hat{t}) \right]  \geq  \text{Var}(\hat{t}).
$$
Substituting terms,

$$
\text{E}\left[\sum_{k \in U} \sum_{l \in \{U: \pi_{kl}>0\}} I_k I_l \frac{\text{Cov}(I_k,I_l)}{\pi_{kl}}\frac{y_k}{\pi_k} \frac{y_l}{\pi_l} +  \sum_{k \in U} \sum_{l \in \{U \backslash k : \pi_{kl}=0\} } \left(I_k \frac{|y_k|^{a_{kl}}}{a_{kl} \pi_k} + I_l \frac{|y_l|^{b_{kl}}}{b_{kl} \pi_l}\right)\right] \geq \text{Var}(\hat{t}).
$$
\hfill $\square$

This estimator is unbiased under a special condition:

\begin{corollary}[Aronow and Samii, 2013, Cor. 1]
If, for all pairs ${k},{l}$ such that $\pi_{kl} = 0$, (i) ${|y_k|^{a_{kl}}}= {|y_l|^{b_{kl}}}$ and (ii) $- y_k y_l = |y_k| |y_l|$,
$$
\text{E} \left[\widehat{\text{Var}_C}(\hat{t}) \right] = \text{Var}(\hat{t}).
$$
\end{corollary}

\begin{proof}
By (i), (ii) and Young's inequality,
$$
\frac{|y_k|^{a_{kl}}}{a_{kl}} + \frac{|y_l|^{b_{kl}}}{b_{kl}} = |y_k| |y_l| = -y_k y_l.
$$
Therefore,
\begin{align}
A^* & =  \sum_{k \in U} \sum_{l \in \{U \backslash k : \pi_{kl}=0\} }  \frac{|y_k|^{a_{kl}}}{a_{kl}} + \frac{|y_l|^{b_{kl}}}{b_{kl}} \nonumber \\ 
& =  \sum_{k \in U} \sum_{l \in \{U \backslash k : \pi_{kl}=0\} } |y_k| |y_l| \nonumber \\ 
& = \sum_{k \in U} \sum_{l \in \{U \backslash k : \pi_{kl}=0\} } -y_k y_l = -A \nonumber
\end{align}

It follows that
$$
 \text{Var}(\hat{t}) + A + A^* = \text{Var}(\hat{t})
$$
and
$$
\text{E} \left[\widehat{\text{Var}_C}(\hat{t}) \right] = \text{Var}(\hat{t}).
$$
\end{proof}

In general, it would be difficult to assign optimal values of $a_{kl}$ and $b_{kl}$ for all pairs ${k},{l}$ such that $\pi_{kl} = 0$. Instead, we examine one intuitive case, assigning all $a_{kl} = b_{kl} = 2$:
$$
\widehat{\text{Var}_{C2}}(\hat{t}) =  \sum_{k \in U} \sum_{l \in \{U: \pi_{kl}>0\}} I_k I_l \frac{\text{Cov}(I_k,I_l)}{\pi_{kl}}\frac{y_k}{\pi_k} \frac{y_l}{\pi_l} +  \sum_{k \in U} \sum_{l \in \{U \backslash k : \pi_{kl}=0\} } \left(I_k \frac{y_k^2}{2\pi_k} + I_l \frac{y_l^2}{2\pi_l}\right).
$$
As a special case of $ \widehat{\text{Var}_{C}}(\hat{t})$, $\widehat{\text{Var}_{C2}}(\hat{t})$ is also conservative:
\begin{corollary}[Aronow and Samii, 2013, Cor. 2]
$$
\text{E} \left[\widehat{\text{Var}_{C2}}(\hat{t})\right] \geq \text{Var}(\hat{t}).
$$
\end{corollary}
\begin{proof}
For all pairs ${k},{l}$ such that $\pi_{kl} = 0$, $\frac{1}{a_{kl}} + \frac{1}{b_{kl}} = \frac{1}{2} + \frac{1}{2} = 1$. Proposition \ref{prop:as-prop2} therefore holds.
\end{proof}

\subsection{Proof of Proposition \ref{prop:tau-consistent}}
We follow the logic of \citet{robinson82}.  $\widehat{\mu_{HT}}(d_k)$ is unbiased for $\mu(d_k)$, and thus we need only consider the variance. 
Condition \ref{cond1} implies that for all values $i$ and $d_k$, $|y_i(d_k)| / \pi_i(d_k) \leq c_3 < \infty$. 
Substituting from Equation \eqref{eq:total_variance}, $N^2 \Var(\widehat{\mu_{HT}}(d_k)) \leq c_3^2  N +  c_3^2 \sum_{i=1}^{N} \sum_{j=1}^{N} g_{ij}.$ 
%
%
Consistency of $\widehat{\mu_{HT}}(d_k)$ for $\mu(d_k)$ is therefore ensured when $\sum_{i=1}^{N} \sum_{j=1}^{N} g_{ij} = o(N^{2})$, as this implies that $\widehat{\mu_{HT}}(d_k) - \mu_{HT}(d_k)  \overset{p}{\longrightarrow} 0$. Consistency of $\widehat{\tau_{HT}}(d_k,d_l)$ for $\tau(d_k,d_l)$ follows by Slutsky's Theorem.

\subsection{Proof of Proposition \ref{prop:tau-interval}}

We follow a proof technique similar to that of \citet{aronow_samii_assenova_dyadic} to establish 
convergence
of $N \widehat{\Var} [ \widehat {\tau_{HT}}(d_k,d_l)]$, though 
for a considerably more general setting.
By Proposition \ref{consvar}, $\E[N \widehat{\Var}[\widehat {\tau_{HT}}(d_k,d_l)]] \geq  N \Var[\widehat {\tau_{HT}}(d_k,d_l)]$. Thus, by Chebyshev's Inequality, $\Var[N \widehat{\Var}[\widehat {\tau_{HT}}(d_k,d_l)]] \arrowp 0$ is sufficient to establish convergence of $N \widehat{\Var}[\widehat {\tau_{HT}}(d_k,d_l)]$ to a value greater than or equal to $N \Var[\widehat {\tau_{HT}}(d_k,d_l)]$, which is itself nonzero by Condition \ref{cond:nonzero-var}. Denote $a_{ij}(D_i,D_j)$ as the sum of the elements in $\widehat{\Var}[\widehat {\tau_{HT}}(d_k,d_l)]$ that incorporate observations $i$ and $j$. Note that all $a_{ij}(D_i,D_j)$ are bounded above by some finite constant by Condition \ref{cond1}. \begin{align*}
\Var[N \widehat{\Var}[\widehat {\tau_{HT}}(d_k,d_l)]] & \leq N^{-2} \Var\left[\sum_{i=1}^{N} \sum_{j=1}^{N} h_{ij} a_{ij}(D_i, D_j)\right] \\
& = N^{-2} \sum_{i=1}^{N} \sum_{j=1}^{N} \sum_{k=1}^{N}  \sum_{l=1}^{N} \Cov[h_{ij} a_{ij}(D_i, D_j), h_{kl} a_{kl}(D_k, D_l)].
\end{align*}
Note that  $\Cov[h_{ij} a_{ij}(D_i, D_j), h_{kl} a_{kl}(D_k, D_l)] \neq 0$ if and only if $h_{ij} = 1$ and $h_{kl} = 1$ and either: 
$h_{ik}  = 1$,  $h_{il} = 1$, $h_{jk}  = 1$, or $h_{jl} = 1$. By Condition \ref{cond:local-dep}, given $m \ll N$, each of these four conditions is satisfied by fewer than $N m^3$ of the elements of the quadruple summation, and the number of elements in their union is at most $4N m^3$. Thus, $\Var[N \widehat{\Var}[\widehat {\tau_{HT}}(d_k,d_l)]] = O(N^{-2}\times N) = O(N^{-1})$. 

Define
$$
t = \frac{\widehat {\tau_{HT}}(d_k,d_l)- {\tau_{HT}}(d_k,d_l)}{\sqrt{\widehat{\Var}[\widehat {\tau_{HT}}(d_k,d_l)]}} = \frac{\widehat {\tau_{HT}}(d_k,d_l)- {\tau_{HT}}(d_k,d_l)}{\sqrt{\Var[\widehat {\tau_{HT}}(d_k,d_l)]}} \left( \frac{\Var[\widehat {\tau_{HT}}(d_k,d_l)]}{\widehat{\Var}[\widehat {\tau_{HT}}(d_k,d_l)]} \right)^{1/2}.
$$
Under Conditions \ref{cond1}, \ref{cond:local-dep}, and \ref{cond:nonzero-var}, then, by \citet[Theorem 2.7]{chenshao2004}, $\frac{\widehat {\tau_{HT}}(d_k,d_l)- {\tau_{HT}}(d_k,d_l)}{\sqrt{\Var[\widehat {\tau_{HT}}(d_k,d_l)]}}$ is asymptotically $\N(0,1)$, while $( \Var[\widehat {\tau_{HT}}(d_k,d_l)]/\widehat{\Var}[\widehat {\tau_{HT}}(d_k,d_l)])^{1/2}$ converges in probability to a quantity in $(0,1]$. By Slutsky's Theorem, $t$ is asymptotically normal and Wald-type confidence intervals constructed as $\widehat {\tau_{HT}}(d_k,d_l) \pm z_{1-\alpha/2}\sqrt{\widehat{\Var}[\widehat {\tau_{HT}}(d_k,d_l)]}$ will tend to cover $\tau_{HT}(d_k,d_l)$ at least $100(1-\alpha)\%$ of the time as $N\rightarrow\infty$.

\subsection{Proof of Proposition \ref{prop:misspec}}
The result follows from iterating expectations:
\begin{align}
\E[\widehat{y^T_{HT}(d_k)}] & =  \E\left[ \sum_{i=1}^N \mathbf{I}(D_i = d_k) \frac{Y_i}{\pi_i(d_k)} \right] \nonumber \\
& = \sum_{i=1}^N \E\left[ \frac{\mathbf{I}(D_i = d_k)}{\pi_i(d_k)} \E\left[ Y_i | D_i = d_k\right]\right] \nonumber \\
& = \sum_{i=1}^N \E\left[ Y_i | D_i = d_k \right] \nonumber \\
& = \sum_{i=1}^N \frac{\sum_{\mathbf{z}:f(\mathbf{z},\theta_i)=d_k} p_{\mathbf{z}}y_i^r(\mathbf{z})}{\sum_{\mathbf{z}:f(\mathbf{z},\theta_i)=d_k} p_{\mathbf{z}}}\nonumber \\
& = \sum_{i=1}^N \sum_{\mathbf{z}:f(\mathbf{z},\theta_i)=d_k} \frac{p_{\mathbf{z}}}{\pi_i(d_k)}y_i^r(\mathbf{z}).\nonumber
\end{align}

\clearpage
\section{Simulation study exposure probabilities}
\begin{figure}[!ht]
\begin{center}
\includegraphics[width=1\textwidth]{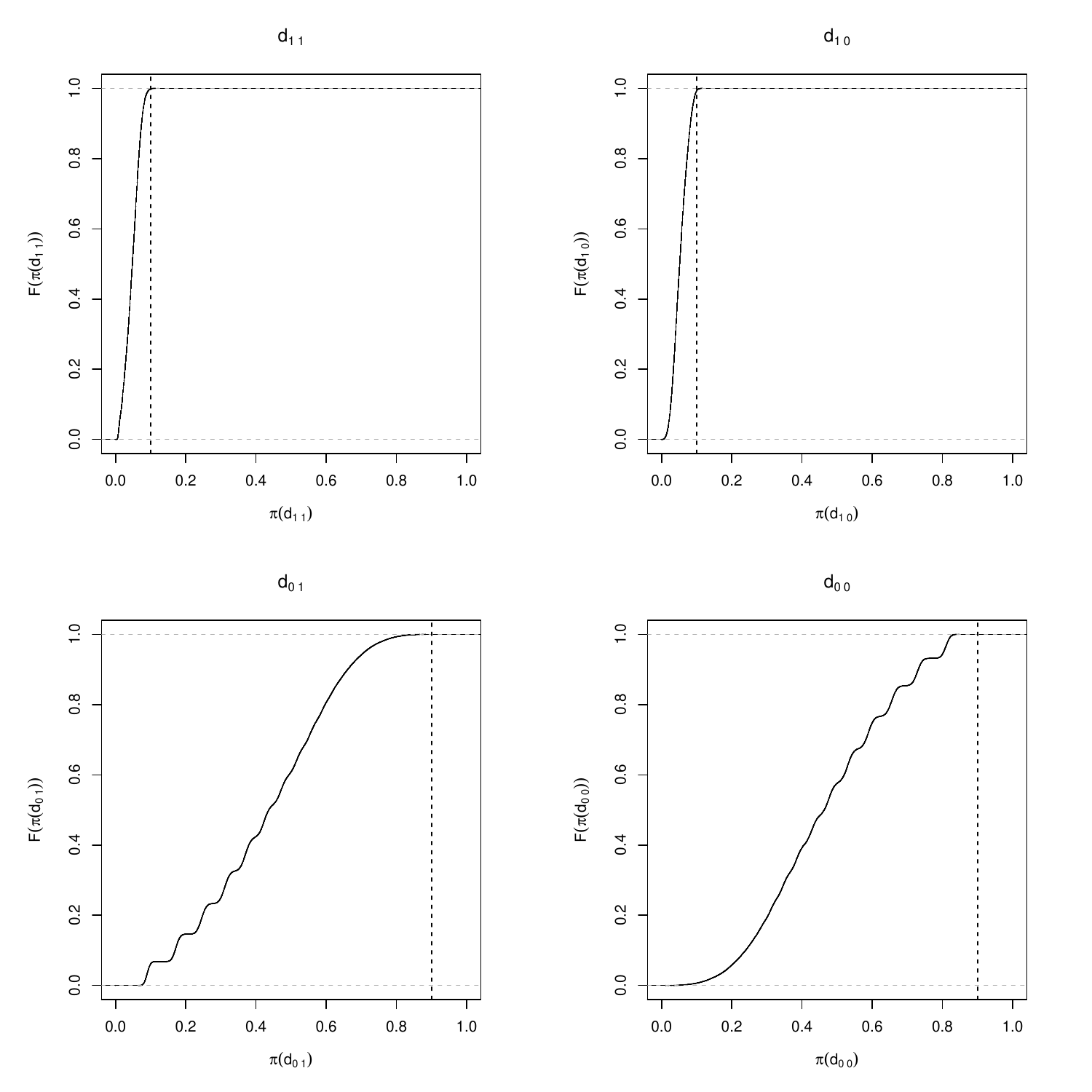}
\caption{Empirical CDFs of probabilities for the four types of exposure in the simulated social network experiment.}
\label{ecdfs}
\end{center}
\end{figure}
\clearpage

\bibliography{interference-aoas-resubmit-2016-07-11.bib}

\end{document}